\documentclass[11pt]{article}
 \usepackage[top=1in,bottom=1in,left=1.5in,right=1.5in]{geometry}
  \usepackage{amsmath,amssymb}
  \usepackage{latexsym}
  \usepackage{setspace}
  \usepackage{tikz}

  \newcommand{\n}{\mathbf{n}}

  \newcommand{\hs}{\hspace*{\parindent}}

  \newcommand{\qed}{\hspace*{\fill} $\Box$\\}

  \newtheorem{theo}{\bfseries \hs Theorem}[section]

  \numberwithin{equation}{section} 

\begin{document}
\title{Upper Bound}
\author{Elliot Krop\thanks{Department of Mathematics, Clayton State University, (ElliotKrop@clayton.edu)} \and Irina Krop\thanks{Depaul University, (irina.krop@gmail.com)}}
\date\today
\maketitle

 We show the improved upper bound for the bipartite problem, when the two parts of $G$ are of equal size.
 
 \begin {theo}
 
 \[f(K_{n,n},G_4,3) \leq n, \mbox{ for all }n\geq 3\]

 \end {theo}

 \subsection {The Coloring}
\vspace{.1 in}
 
 We will explore the matrix
 
 $$G=\begin{pmatrix}
1           &2        &3         &. &r         &. &c+1        &.     &n  \cr
3           &1        &2         &. &r-1       &. &c          &.     &n-1 \cr
v_{3}     &n-1      &1         &. &r-2       &. &c-1        &.     &n-2 \cr
.           &.        &.         &. &.         &. &.          &.     &.   \cr
v_{n+1-r} &r+1      &r+2       &. &1         &. &r+c        &.     &r    \cr
.           &.        &.         &. &.         &. &.          &.     &.    \cr
v_{n-1}   &3        &4         &. &r+1       &. &c+2        &.     &2    \cr
n-2         &u_{2}  &u_{3}   &. &u_{r} &. &u_{c+1}  &.     &1
\end{pmatrix}$$

The values of $v_i$ and $u_i$ will be defined shortly.

\vspace{.1 in}
 
Let permutation $\sigma$ be the $n-1$ cycle $(1\; 2\; \ldots\;  n-1)$. That is, $\sigma$ sends $i$ to $i+1 \pmod{n-1}$. For a natural number $m$ we shall write $m \pmod {n-1}$ for its representative in $\{1,2,\ldots ,n-1\}$. For each $r$ we define $\sigma^{(r)}$ by the rule $\sigma^{(r)}{(c)} \equiv r+c \pmod{n-1}$.  Let us start with the matrix

$$C=\begin{pmatrix}
2                   &3                      &. &c+1                 &.       &n                   \cr
\sigma^{0}{(1)}   &\sigma^{0}{(2)}      &. &\sigma^{0}{(c)}   &.       &\sigma^{0}{(n-1)} \cr
\sigma^{n-2}{(1)} &\sigma^{n-2}{(2)}    &. &\sigma^{n-2}{(c)} &.       &\sigma^{n-2}{(n-1)}\cr
.                   &.                      &. &.                   &.       &. \cr
\sigma^{r}{(1)}   &\sigma^{r}{(2)}      &. &\sigma^{r}{(c)}   &.       &\sigma^{r}{(n-1)}\cr
.                   &.                      &. &.                   &.       &.  \cr
\sigma^{2}{(1)}   &\sigma^{2}{(2)}      &.   &\sigma^{2}{(c)} &.      &\sigma^{2}{(n-1)}
\end{pmatrix}$$

\vspace{.1 in}

We define matrix $G$ by adding the first column $V=\{v_1,\dots,v_{n-1},vu\}$ and the last row $U=\{vu,u_2,\dots,u_n\}$ to the matrix $C$.

$$G=\begin{pmatrix}
v_{1}       &2                   &3                      &. &c+1                 &.       &n                   \cr
v_{2}       &\sigma^{0}{(1)}     &\sigma^{0}{(2)}        &. &\sigma^{0}{(c)}     &.       &\sigma^{0}{(n-1)} \cr
v_{3}       &\sigma^{n-2}{(1)}   &\sigma^{n-2}{(2)}      &. &\sigma^{n-2}{(c)}   &.       &\sigma^{n-2}{(n-1)}\cr
.             &.                   &.                      &. &.                   &.       &. \cr
v_{n+1-r}   &\sigma^{r}{(1)}     &\sigma^{r}{(2)}        &. &\sigma^{r}{(c)}     &.       &\sigma^{r}{(n-1)}\cr
.             &.                   &.                      &. &.                   &.       &.  \cr
v_{n-1}     &\sigma^{2}{(1)}    &\sigma^{2}{(2)}         &. &\sigma^{2}{(c)}     &.       &\sigma^{2}{(n-1)}\cr
vu            &u_{2}                &u_{3}                  &. &u_{c+1}           &.       &u_{n}
\end{pmatrix}$$

\vspace{.1 in}

The entries of $G$ will be defined as follows:
for every 4-tuple $(i,j;l,m)$ with $1 \leq i<j \leq n$ and $1\leq l<m \leq n$ the  ($2\times 2$) matrix

$$G(i,j;l,m)=\begin{pmatrix}
                    a_{il}&a_{im}\cr
                    a_{jl}&a_{jm}\end{pmatrix}$$

\vspace{.1 in}
 
 We consider the colorings for the edges $V$ and $U$ in three types of even $n\pmod 6$.

\vspace{.1 in}
 
 \noindent \emph{\textbf{Type 1:} Matrix $G_1=G$ for $n \equiv 2 \pmod 6$};    $ \;  \; \left[ \;  n=2+6k, \; k\geq 1 \; \right]$

   \[ a_{i,1} = \left\{\begin{array}{ll}
    1, & i=1\\
    3, & i=2\\
    n, & 3\leq i \leq \frac{n}{2} +1\\
    2(i-1)-n, & \frac{n}{2}+2\leq i\leq n-1\\ 
    n-2, & i=n  
    \end{array}\right. \]

    \[ a_{n,l} = \left\{\begin{array}{ll}
    n-2l, & 1\leq l \leq \frac{n}{2}-1\\
    n, & \frac{n}{2} \leq l \leq n-2\\
    n-1, & l=n-1\\
    1, & l=n    
    \end{array}\right. \]
 
     \vspace{.1 in}
     
\noindent \emph{\textbf{Type 2:} Matrix $G_2=G$ for $n \equiv 6 \pmod 6$};  $ \; \; \left[ \; n=6+6k, \; k\geq 1 \; \right]$

 \noindent  We define $Y$ as $\frac{n}{2}-2$ for even $k$,  and as $\frac{n}{2}+1$ for odd $k$.
   
   \[ a_{i,1} = \left\{\begin{array}{ll}
    1, & i=1\\
    3, & i=2\\
    n, & 3\leq i \leq \frac{n}{2} +1\\
    Y, & i=\frac{n}{2}+2\\
    2(i-2)-n, & \frac{n}{2}+3\leq i\leq n-1\\
    n-2, & i=n   
    \end{array}\right. \]

   \[ a_{n,l} = \left\{\begin{array}{ll}
    n-2, & l=1\\
    n-2(l+1), & 2\leq l \leq \frac{n}{2}-2\\
    Y, & l= \frac{n}{2}-1\\
    n, & \frac{n}{2} \leq l \leq n-2\\
    n-1, & l=n-1\\
    1, & l=n    
    \end{array}\right. \]
 \vspace{.1 in}
 
 \noindent  Exception for $n=6$;  $\left[ k= 0 \right]$ the first row $V=\{1,5,6,6,4,\}$, the last column $U=\{3,6,6,6,5,1\}$.   
    
   \vspace{.1 in}

 \noindent \emph{\textbf{Type 3:} Matrix $G_3=G$ for $n \equiv 4 \pmod 6$};  $ \; \; \left[ \; n=4+6k, \; k\geq 4 \; \right]$

  \noindent  The regularity  starts with $n>22$.  
  
      \[ a_{i,1} = \left\{\begin{array}{ll}
    1, & i=1\\
    3, & i=2\\
    n, & 3\leq i \leq \frac{n}{2} +1\\
    n-9, & i=\frac{n}{2}+2\\
    2(i-2)-n, & \frac{n}{2}+3\leq i\leq \frac{5n+4}{6}\\
    2(i-1)-n, & \frac{5n+10}{6}\leq i\leq n-1\\
    n-2, & i=n   
    \end{array}\right. \]

     \[ a_{n,l} = \left\{\begin{array}{ll}
    n-2l, & 1\leq l \leq \frac{n-4}{6}\\
    n-2(l+1), & \frac{n+2}{6}\leq l \leq \frac{n}{2}-2\\
    n-9, & l=\frac{n}{2}-1\\
    n, & \frac{n}{2} \leq l \leq n-2\\
    n-1, & l=n-1\\
    1, & l=n    
    \end{array}\right. \]

 \vspace{.1 in}
 
Exceptions:

 For $n=10$ we replace $(n-9)$ with $(n-8)$. 

 For $n=16$ we replace $(n-9)$ with $(n-11)$.
  
 For $n=22$ we replace 
 $(n-9)$ with $(n-5)$ and the definitions:
 
  \[ a_{i,1} = \left\{\begin{array}{ll}
    2(i-2)-n, & \frac{n}{2}+3\leq i\leq \frac{5n-2}{6}\\
    2(i-1)-n, & \frac{5n+4}{6}\leq i\leq n-1
    \end{array}\right. \]

    \[a_{n,l} =  \left\{\begin{array}{ll}   
     n-2l, & 1\leq l \leq \frac{n-10}{6}\\
     n-2(l+1), & \frac{n-4}{6}\leq l \leq \frac{n}{2}-2
    \end{array}\right. \]

 \vspace{.1 in}
 
 \subsection{Proof in Steps}
 
 \begin{proof}
 
  First, we show that every $4-cycle$ defined in the \emph{basic coloring} is $almost$ $rainbow$.  That is, given $i<j$ and $l<m$ we show that $a_{i,l},a_{j,l},a_{i,m},a_{j,m}$ contains at least three distinct elements in the \emph{basic coloring}.
  
\vspace{.2 in}

\noindent \emph{\textbf{Step 1}}:

We start with the matrix $C$ and look at two occurrences, which are identical for each of the types of even $n\pmod 6$ specified above.
\vspace{.1 in}

\emph{\textbf{Case 1}}:
We take the submatrix of $G(i,j;l,m)$ with $2 \leq l < m \leq n$, $2 \leq i <j < n$, and let $s=(n+1)-i$, $t=(n+1)-j$. 
A typical ($2\times 2$) submatrix in this case has the form:

$$\begin{pmatrix}
                    \sigma^{s}{(l-1)}& \sigma^{s}{(m-1)}\cr
                    \sigma^{t}{(l-1)}& \sigma^{t}{(m-1)}\end{pmatrix}$$

We wish to show there are three distinct elements:
$\sigma^{s}{(l-1)}\neq \sigma^{t}{(l-1)}$,
$\sigma^{t}{(l-1)}\neq \sigma^{t}{(m-1)}$, 
$\sigma^{s}{(l-1)}\neq \sigma^{t}{(m-1)}$.
Suppose $\sigma^{(s)}{(l-1)} \equiv \sigma^{(t)}{(l-1)}  \Rightarrow s\equiv t$ $\pmod {n-1}$, which is a contradiction.
Suppose $\sigma^{(t)}{(l-1)} \equiv \sigma^{(t)}{(m-1)} \Rightarrow l\equiv m$ $\pmod {n-1}$, which is a contradiction.
Suppose $\sigma^{(s)}{(l-1)} \equiv \sigma^{(t)}{(m-1)} \Rightarrow s+l\equiv t+m$ $\pmod {n-1}$, and assume there are three distinct elements:
$\sigma^{t}{(l-1)}\neq \sigma^{t}{(m-1)}$,
$\sigma^{s}{(m-1)}\neq \sigma^{t}{(m-1)}$, 
$\sigma^{s}{(m-1)}\neq \sigma^{t}{(l-1)}$. 
Follow the argument above the first two inequalities are correct.
Suppose $\sigma^{s}{(m-1)} \equiv \sigma^{t}{(l-1)}\Rightarrow s+m\equiv t+l$ $\pmod {n-1}$.
Subtracting equations $s+l \equiv t+m$ and  $s+m\equiv t+l\Rightarrow l \equiv m$ $\pmod {n-1}$, which is a contradiction.
One of the following two sets has three distinct elements:
$\{ \sigma^{s}{(l-1)},\;   \sigma^{t}{(l-1)},\;   \sigma^{t}{(m-1)} \}$ or $\{ \sigma^{t}{(l-1)},\;   \sigma^{t}{(m-1)},\;   \sigma^{s}{(m-1)} \}$. 
  
\vspace{.1 in}

\emph{\textbf{Case 2}}:
We take the submatrix of $G(i,j;l,m)$ with $2 \leq l< m \leq n$, $i=1$, $1<j<n$, and let $r=(n+1)-j$. 
A typical ($2\times 2$) submatrix has the form:

$$\begin{pmatrix}
                    l&m\cr
                    \sigma^{r}{(l-1)}& \sigma^{r}{(m-1)}\end{pmatrix}$$

We wish to show there are three distinct elements: 
$l \neq m$,
$m \neq \sigma^{r}{(m-1})$, 
$l \neq \sigma^{r}{(m-1)}$.  
Suppose $m \equiv \sigma^{(r)}{(m-1)}\Rightarrow r \equiv 1$ $\pmod {n-1}$, which is a contradiction.
Suppose $l \equiv \sigma^{(r)}{(m-1)}\Rightarrow l \equiv r+m-1$ $\pmod {n-1}$, and assume there are three distinct elements: 
$\sigma^{r}{(l-1)}\neq \sigma^{r}{(m-1)}$,
$m\neq \sigma^{r}{(m-1)}$,
$m\neq \sigma^{r}{(l-1)}$.  
The first two inequalities are correct. Suppose $m \equiv \sigma^{r}{(l-1)}\Rightarrow m \equiv r+l-1$ $\pmod {n-1}$. 
Subtracting equations $l \equiv r+m-1$ and  $m \equiv r+l-1\Rightarrow r=1$ $\pmod {n-1}$, which is a contradiction.
One of the following two sets has three distinct elements:
$\{l,\;  m,\;  \sigma^{r}{(m-1)} \}$ or $\{ \sigma^{r}{(l-1)},\; \sigma^{r}{(m-1)},\;  m \}$. 
 
\vspace{.2 in}
 
\noindent \emph{\textbf{Step 2}}:

For matrix $G(i,j;l,m)$ with  $i=1$, $j=n$ and $2 \leq l <m\leq n$ we look at five cases and consider every matrix type defined above of even $n\pmod 6$.

\vspace{.1 in}

\emph{\textbf{Case 1}}: We take $G(i,j;l,m)$ with   $2\leq l<m \leq \frac{n}{2}-1$, $i=1$, $j=n$. 

\vspace{.1in}

\noindent \textit{Consider  $G_1$}.

  $$\begin{pmatrix}
       l&m\cr
       n-2l&n-2m\end{pmatrix}$$
                   
We wish to show there are three distinct entries:
$l\neq n-2l$,
$n-2l\neq n-2m$, 
$l\neq n-2m$. 
Suppose $l =n-2l\Rightarrow 3l=n$ and since $n=2+6k$  this is a contradiction.
Suppose $n-2l=n-2m \Rightarrow l=m$, which is a contradiction.
Suppose $l=n-2m$ and we wish to show there are three distinct elements:
$l\neq m$,
$l\neq n-2l$,
$m\neq n-2l$.
As shown above the first inequality is correct.
Suppose $m=n-2l$ and since $l=n-2m \Rightarrow l=m$, which is a contradiction.
One of the following two sets has three distinct elements: $\{l, n-2l, n-2m\}$ or $\{l, m, n-2l\}$. 

\vspace{.1in}

\noindent \textit{Consider $G_2$}

{\bf 1.}  $G_2$ with  $2 \leq l< m \leq \frac{n}{2}-2$, $i=1$, $j=n$.

$$\begin{pmatrix}
                   l&m\cr
                   n-2(l+1)&n-2(m+1)\end{pmatrix}$$

We wish to show there are three distinct entries:
$l\neq n-2(m+1)$,
$l\neq n-2(l+1)$, 
$n-2(l+1)\neq n-2(m+1)$.
Suppose $l =n-2l-2\Rightarrow 3l=n-2$ and since $n=6+6k$  this is a contradiction.
Suppose $n-2l=n-2m \Rightarrow l=m$, which is a contradiction.
Suppose $l=n-2(m+1)$ and we wish to show there are three distinct elements:
$l\neq m$,
$l\neq n-2(l+1)$,
$m\neq n-2(l+1)$.
As shown above the first inequality is correct.
Suppose $m=n-2(l+1)$ and since $l=n-2(m+1) \Rightarrow l=m$, which is a contradiction.
One of the following two sets has three distinct elements: $\{l, n-2(m+1), n-2(l+1)\}$ or $\{l, m, n-2(l+1)\}$. 

{\bf 2.} $G_2$ with $ 2 \leq l \leq \frac{n}{2}-2$, $m=\frac{n}{2}-1$, $i=1$, $j=n$.

$$\begin{pmatrix}
                   l&m\cr
                   n-2(l+1)&Y\end{pmatrix}$$

If  \textit{$K$ is Even}  $\Rightarrow m=\frac{n}{2}-1$, $Y=\frac{n}{2}-2 \Rightarrow Y=m-1$.

We wish to show there are three distinct entries:
$l\neq m$,
$m\neq m-1$, 
$l\neq m-1$.
Assume $l= m-1$ and we wish to show there are three distinct entries:
$m\neq m-1$,
$m\neq n-2(l+1)$,
$m-1\neq n-2(l+1)$.
Suppose $m= n-2(l+1)$ and since $l=m-1$ and $m=\frac{n}{2}-1 \Rightarrow  n=6$, which is a contradiction.
Suppose $m-1= n-2(l+1)$ and since $m=\frac{n}{2}-1$ and $l= m-1 \Rightarrow  n=8$, which is a contradiction.
One of the following two sets has three distinct elements: $\{l, m, Y\}$ or $\{m, Y, n-2(l+1)\}$. 
\vspace{.1in}

If  \textit{$K$ is Odd} $ \Rightarrow m=\frac{n}{2}-1$, $Y=\frac{n}{2}+1 \Rightarrow  Y=m+2$.
Three distinct elements are $\{l, m, Y\}$. 

\vspace{.1in}

\noindent \textit{Consider $G_3$}

{\bf 1.}  $G_3$ with $i=1$, $j=n$ and ($2 \leq l<m \leq \frac{n-4}{6}$ or $\frac{n+2}{6} \leq l<m \leq \frac{n}{2}-2$). The argument is similar to above one with $G_1$. One of the following two sets has three distinct elements: $\{l, n-2l, n-2m \}$ or $\{ l, m, n-2l \}$. 

{\bf2.} $G_3$ with   $2 \leq l \leq \frac{n-4}{6}$, $\frac{n+2}{6} \leq m \leq \frac{n}{2}-2$, $i=1$, $j=n$.

$$\begin{pmatrix}
                   l&m\cr
                   n-2l&n-2(m+1)\end{pmatrix}$$

We wish to show there are three distinct entries:
$l\neq n-2l$,
$n-2l\neq n-2(m+1)$, 
$l\neq n-2(m+1)$.
Suppose $l =n-2l\Rightarrow 3l=n$ and since $n=4+6k$ this is a contradiction.
Suppose $n-2l= n-2(m+1)\Rightarrow l=m+1$, which is a contradiction.
Suppose $l=n-2(m+1)$ and we wish to show there are three distinct elements:
$l \neq m$,
$l\neq n-2l$,
$m \neq n-2l$.
The first two inequalities are correct.
Suppose $m = n-2l$ and since  $l=n-2(m+1)  \Rightarrow m=l-2$, which is a contradiction.
One of the following two sets has three distinct elements: $\{l, n-2l, n-2(m+1))\}$ or $\{ l, m, n-2l \}$. 

{\bf3.}  $G_3$ with $2 \leq l \leq \frac{n-4}{6}$, $m=\frac{n}{2}-1$, $i=1$, $j=n$.

 $$\begin{pmatrix}
                   l&m\cr
                   n-2l&n-9\end{pmatrix}$$

We wish to show  there are three distinct entries:
$l\neq m$,
$m\neq n-9$,
$l\neq n-9$.
Suppose $l= n-9$ and since  $l <\frac{n-4}{6}\Rightarrow n-9<\frac{n-4}{6}\Rightarrow n<10$, which is a contradiction.
Suppose $m =n-9\Rightarrow \frac{n}{2}-1=n-9\Rightarrow n=16$, which is a contradiction.
There are three distinct elements $\{l, m, n-9\}$. 

{\bf  4.}  $G_3$ with  $\frac{n+2}{6}\leq l \leq \frac{n}{2}-2$, $m=\frac{n}{2}-1$, $i=1$, $j=n$.

 $$\begin{pmatrix}
                   l&m\cr
                   n-2(l+1)&n-9\end{pmatrix}$$

We wish to show  there are three distinct entries:
$l\neq m$,
$m\neq n-9$,
$l\neq n-9$.
Suppose $l= n-9$ and since  $l <\frac{n}{2}-2\Rightarrow n-9<\frac{n}{2}-2\Rightarrow n<14$, which is a contradiction.
Suppose $m =n-9\Rightarrow \frac{n}{2}-1=n-9\Rightarrow n=16$, which is a contradiction.
There are three distinct elements $\{l, m, n-9\}$. 

\vspace{.1in}

Exceptions:
For $n=16$ three distinct elements  are $\{14-2l, 5, 7\}$.
For $n=22$ three distinct elements are $\{l, m, 17\}$

\vspace{.1in}

\emph{\textbf{Case 2}}:  For matrix $G(i,j;l,m)$ with $i=1$, $j=n$ and ($\frac{n}{2} \leq l < m \leq n-2$ or $2\leq l \leq \frac{n}{2}-1$, $\frac{n}{2} \leq  m \leq n-2$) three distinct elements are $\{l, m, n\}$.

\vspace{.1in}

\emph{\textbf{Case 3}}: We take the submatrix $G(i,j;l,m)$ with  $2\leq l \leq \frac{n}{2}-1$, $ m = n-1$, $i=1$, $j=n$.
\vspace{.1in}

\noindent \textit{Consider $G_1$}

$$\begin{pmatrix}
                   l&m\cr
                   n-2l&n-1\end{pmatrix}$$

We wish to show there are three distinct entries:
$l \neq n-1$,
$n-2l \neq n-1$,
$l \neq n-2l$.
Suppose  $n-2l = n-1\Rightarrow l=\frac{1}{2}$, which is a contradiction.
Suppose  $l=n-2l\Rightarrow 3l=n$ and since $n=2+6k$ this is a contradiction.
There are three distinct elements $\{l, n-1, n-2l\}$.
\vspace{.1in}

\noindent \textit{Consider $G_2$}

{\bf 1.} $G_2$ with $2\leq l \leq \frac{n}{2}-2$, $m =n-1$, $i=1$, $j=n$.

$$\begin{pmatrix}
                   l&m\cr
                   n-2(l+1)&n-1\end{pmatrix}$$

We wish to show there are three distinct entries:                            
$l \neq n-1$,
$n-2(l+1) \neq n-1$,
$l \neq n-2(l+1)$.
Suppose  $n-2(l+1) = n-1 \Rightarrow l=-\frac{1}{2}$, which is a contradiction.
Suppose $l=n-2(l+1)\Rightarrow 3l=n-2$ and since $n=6+6k$   this is a contradiction.
There are three distinct elements $\{l, n-1, n-2(l+1)\}$.

{\bf 2.}  $G_2$ with  $l =\frac{n}{2}-1$, $m=n-1$, $i=1$, $j=n$.

$$\begin{pmatrix}
                   l&m\cr
                   Y&n-1\end{pmatrix}$$

If  \textit{$K$ is Even} $\Rightarrow Y=\frac{n}{2}-2$.                         
If  $n-1 = Y\Rightarrow n-1=\frac{n}{2}-2\Rightarrow n=-2$, which is a contradiction. Three distinct entries are $\{l,Y,n-1\}$.
If  \textit{$K$ is Odd} $\Rightarrow Y=\frac{n}{2}+1$, and three distinct entries are $\{l,Y,n-1\}$.

\vspace{.1in}

\noindent \textit{Consider $G_3$}

{\bf 1.} For $G_3$ with $l \leq \frac{n-4}{6}$, $ m = n-1$, $i=1$, $j=n$ three distinct elements are $\{l, n-1, n-2l\}$ (similar to $G_1$.)

{\bf 2.} For $G_3$ with $\frac{n+2}{6} \leq l < n-1$, $m = n-1$, $i=1$, $j=n$.

$$\begin{pmatrix}
                   l&m\cr
                   n-2(l+1)&n-1\end{pmatrix}$$

We wish to show there are three distinct entries:        
$l \neq n-1$,
$n-2(l+1)\neq n-1$,
$l\neq  n-2(l+1)$.
Suppose  $n-2(l+1)= n-1\Rightarrow l=-\frac{1}{2}$, which is a contradiction.
Suppose $l= n-2(l+1)\Rightarrow 3l=n-2$ and since $n=4+6k$ this is a contradiction.
Three distinct elements are $\{l, n-1, n-2(l+1)\}$. 

{\bf 3.} For $G_3$ with $l=\frac{n}{2}-1$, $m=n-1$, $i=1$, $j=n$ three distinct entries are $\{l, n-9, n-1)\}$.

Exceptions: For $n=16$ three distinct entries are $\{7, 5, 15 \}$.
For $n=22$ three distinct entries  are $\{10, 17, 21\}$.
\vspace{.1in}

\emph{\textbf{Case 4}}: For the submatrix $G(i,j;l,m)$ with   $\frac{n}{2} \leq l \leq  n-2$, $ m=n-1$, $i=1$, $j=n$
three distinct elements are $\{l,n,n-1\}$.

\vspace{.1in}

\emph{\textbf{Case 5}}: For the submatrix $G(i,j;l,m)$ with  $l = n-1$, $ m=n$, $i=1$, $j=n$ three distinct elements are $\{n-1,n,1\}$.

\vspace{.2in}

\noindent \emph{\textbf{Step 3}}:

For matrix $G(i,j;l,m)$ with $2\leq i \leq n-1$, $j=n$, $2\leq l <m \leq n$ we look at five cases and consider every matrix type defined above of even $n\pmod 6$.

\vspace{.1 in}

\emph{\textbf{Case 1}}: We take $G(i,j;l,m)$.

\noindent \textit{Consider $G_1$} with $2 \leq l<m \leq \frac{n}{2}-1$, $2\leq i \leq n-1 $, $j=n$.

$$\begin{pmatrix}
                     \sigma^{s}{(l-1)}& \sigma^{s}{(m-1)}\cr
                   n-2l&n-2m\end{pmatrix}$$

We wish to show there are three distinct entries:       
$ \sigma^{s}{(l-1)}\neq  n-2l$,
$\sigma^{s}{(l-1)}\neq  n-2m$,
$ n-2l\neq  n-2m$.
Suppose $ n-2l=n-2m \Rightarrow m =l$, which is a contradiction.
Suppose $ \sigma^{s}{(l-1)}\equiv  n-2l$ and we wish to show there are three distinct entries:       
$n-2m \neq  n-2l$,
$ \sigma^{s}{(m-1)}\neq n-2l$,
$\sigma^{s}{(m-1)}\neq  n-2m$.
The first statement is correct.
Since $\sigma^{s}{(l-1)}\equiv n-2l$ and $ \sigma^{s}{(l-1)}\neq  \sigma^{s}{(m-1)} \Rightarrow  \sigma^{s}{(m-1)}\neq n-2l$.
Subtracting equations $\sigma^{s}{(m-1)}\equiv  n-2m$  and  $ \sigma^{s}{(l-1)}\equiv  n-2l \Rightarrow m \equiv l$ $\pmod {n-1}$, which is a contradiction.
Suppose $ \sigma^{s}{(l-1)}\equiv  n-2m$ and we wish to show there are three distinct entries:  
$n-2m\neq  n-2l$,
$ \sigma^{s}{(m-1)}\neq n-2m$,
$\sigma^{s}{(m-1)}\neq   n-2l$.
The statement $n-2m\neq  n-2l$ is correct.
Since  $\sigma^{s}{(l-1)}\equiv  n-2m$ and $ \sigma^{s}{(l-1)}\neq  \sigma^{s}{(m-1)}\Rightarrow  \sigma^{s}{(m-1)}\neq n-2m$.
Subtracting equations $\sigma^{s}{(m-1)}\equiv  n-2l$  and  $ \sigma^{s}{(l-1)}\equiv  n-2m \Rightarrow m \equiv l$ $\pmod {n-1}$, which is a contradiction.
One of the following two sets has three distinct elements:  
$\{ \sigma^{s}{(l-1)},n-2l, n-2m \}$ or $\{\sigma^{s}{(m-1)}, n-2l, n-2m \}$.
\vspace{.1in}

\noindent \textit{Consider $G_2$}

{\bf 1.} $G_2$ with $2 \leq l <m \leq \frac{n}{2}-2$, $2\leq i \leq n-1 $, $j=n$.
$$\begin{pmatrix}
                    \sigma^{s}{(l-1)}& \sigma^{s}{(m-1)}\cr
                   n-2(l+1)&n-2(m+1)\end{pmatrix}$$
We wish to show there are three distinct entries:    
$ \sigma^{s}{(l-1)}\neq  n-2(l+1)$,
$\sigma^{s}{(l-1)}\neq  n-2(m+1)$,
$ n-2(l+1)\neq  n-2(m+1)$.
Suppose $ n-2(l+1)= n-2(m+1)\Rightarrow m =l$, which is a contradiction.
Suppose $ \sigma^{s}{(l-1)}\equiv  n-2(l+1)$ and we wish to show there are three distinct entries:  
$n-2(m+1) \neq  n-2(l+1)$,
$ \sigma^{s}{(m-1)}\neq n-2(l+1)$,
$\sigma^{s}{(m-1)}\neq n-2(m+1)$.
The first statement is correct. 
Since  $\sigma^{s}{(l-1)}\equiv  n-2(l+1)$ and $ \sigma^{s}{(l-1)}\neq  \sigma^{s}{(m-1)} \Rightarrow  \sigma^{s}{(m-1)}\neq n-2(l+1)$.
Subtracting equations $\sigma^{s}{(m-1)}\equiv  n-2(m+1)$  and $ \sigma^{s}{(l-1)}\equiv  n-2(l+1) \Rightarrow m \equiv l$ $\pmod {n-1}$, which is a contradiction.

Assume $ \sigma^{s}{(l-1)}\equiv  n-2(m+1)$ and we wish to show there are three distinct entries:  
$n-2(m+1)\neq  n-2(l+1)$,
$\sigma^{s}{(m-1)}\neq n-2(m+1)$,
$\sigma^{s}{(m-1)}\neq  n-2(l+1)$
The first statement is correct. 
Since  $\sigma^{s}{(l-1)}\equiv  n-2(m+1)$ and $ \sigma^{s}{(l-1)}\neq  \sigma^{s}{(m-1)}\Rightarrow  \sigma^{s}{(m-1)}\neq n-2(m+1)$.
Subtracting equations  $\sigma^{s}{(m-1)}\equiv  n-2(l+1)$  and $ \sigma^{s}{(l-1)}\equiv  n-2(m+1)\Rightarrow m \equiv l$ $\pmod {n-1}$, which is a contradiction.
One of the following two sets has three distinct elements:  
$\{ \sigma^{s}{(l-1)},n-2(l+1), n-2(m+1)\}$ or $\{ \sigma^{s}{(m-1)}, n-2(l+1), n-2(m+1)\}$.

{\bf  2.}  $G_2$ with $2 \leq l \leq  \frac{n}{2}-2$, $m=\frac{n}{2}-1$, $2\leq i \leq n-1 $, $j=n$.

$$\begin{pmatrix}
                   \sigma^{s}{(l-1)}& \sigma^{s}{(m-1)}\cr
                   n-2(l+1)&Y\end{pmatrix}$$

If  \textit{$K$ is Even} $ \Rightarrow Y=\frac{n}{2}-2$. 
We wish to show there are three distinct entries:  
$ \sigma^{s}{(l-1)} \neq  \sigma^{s}{(\frac{n}{2}-2)}$,
$ \sigma^{s}{(\frac{n}{2}-2)}\neq \frac{n}{2}-2$,
$\sigma^{s}{(l-1)}\neq \frac{n}{2}-2$.
Two first inequalities are correct. 
Assume $\sigma^{s}{(l-1)}\equiv \frac{n}{2}-2$ and we wish to show there are three distinct entries: 
$\sigma^{s}{(\frac{n}{2}-2)} \neq (\frac{n}{2}-2)$,
$ n-2(l+1)\neq \frac{n}{2}-2$,
$\sigma^{s}{(\frac{n}{2}-2)} \neq  n-2(l+1)$.
Suppose  $n-2(l+1)= \frac{n}{2}-2\Rightarrow n=4l$.  Hence $n=6(k+1)\Rightarrow 4l=6(k+1)\Rightarrow l=1.5(k+1)$, and since $k$ is even $\Rightarrow l$ is not integer, which is a contradiction. 
Subtracting equations  $\sigma^{s}{(\frac{n}{2}-2)} \equiv n-2(l+1)$ and $\sigma^{s}{(l-1)}\equiv \frac{n}{2}-2\Rightarrow l \equiv 1$, which is a contradiction.
One of the following two sets has three distinct elements: 
$\{ \sigma^{s}{(l-1)},\sigma^{s}{(m-1)}, \frac{n}{2}-2\}$ or $\{ \sigma^{s}{(m-1)},n-2(l+1), \frac{n}{2}-2 \}$.
\vspace{.1in}

If  \textit{$K$ is Odd} $\Rightarrow Y=\frac{n}{2}+1$. 
We wish to show there are three distinct entries:  
$\sigma^{s}{(m-1)} \neq  Y$,
$\sigma^{s}{(m-1)}\neq  n-2(l+1)$,
$ n-2(l+1)\neq  Y$.
Suppose $ n-2(l+1)= Y \Rightarrow 4l=n-6$. Since $n=6k+6\Rightarrow l=1.5k$ and since $k$ is odd $l$ is not integer, this is a contradiction.
Suppose $\sigma^{s}{(m-1)} \equiv Y \Rightarrow s+\frac{n}{2}-2 \equiv  \frac{n}{2}+1\Rightarrow s\equiv 3$ $\pmod {n-1}$.  
Suppose $\sigma^{s}{(m-1)} \equiv n-2(l+1)\Rightarrow s \equiv  \frac{n}{2}-2l$ $\pmod {n-1}$.
We wish to show there are three distinct entries:  
$\sigma^{s}{(l-1)} \neq  n-2(l+1)$,
$n-2(l+1)\neq \frac{n}{2}+1$,
$\sigma^{s}{(l-1)}\neq \frac{n}{2}+1$.
We already know that $n-2(l+1)\neq \frac{n}{2}+1$.
Suppose $\sigma^{s}{(l-1)} \equiv  \frac{n}{2}+1 \Rightarrow 2s \equiv n+4-2l$ $\pmod {n-1}$. 
Suppose $\sigma^{s}{(l-1)} \equiv  n-2(l+1)\Rightarrow 3l+s \equiv n-1$ $\pmod {n-1}$.
Since  $s\equiv 3$ from the equation $2s \equiv n+4-2l \Rightarrow l\equiv \frac{n}{2}-1$, which  is a contradiction.
And from the  equation $s+3l \equiv n-1$ since $n=6k+6 \Rightarrow l\equiv 2k+\frac{2}{3}$, which is a contradiction.
Going back to the assumption $ 2s\equiv n-4l$  and subtracting it from $2s \equiv n+4-2l \Rightarrow l\equiv -2$, which  is a contradiction.
Subtracting $ 2s\equiv n-4l$  from $s+3l \equiv n-1 \Rightarrow l\equiv \frac{n}{2}-1$, which  is a contradiction.
One of the following two sets has three distinct elements: 
$\{Y, \sigma^{s}{(l-1)}, n-2(l+1)\}$ or $\{\sigma^{s}{(m-1)},n-2(l+1), Y\}$.

\vspace{.1in}

\noindent \textit{Consider $G_3$}.

{\bf  1.}  $G_3$ with ($2 \leq l< m \leq \frac{n-4}{6}$ or $\frac{n+2}{6} \leq l<m \leq \frac{n}{2}-2$), $2\leq i \leq n-1 $, $j=n$. This case is similar to the case for the matrix $G_1$.

{\bf 2.}  $G_3$ with  $2\leq l \leq \frac{n-4}{6}$, $\frac{n+2}{6}\leq m \leq \frac{n}{2}-2$, $2\leq i \leq n-1 $, $j=n$.

$$\begin{pmatrix}
                   \sigma^{s}{(l-1)}&\sigma^{s}{(m-1)}\cr
                   n-2l&n-2(m+1)\end{pmatrix}$$

We wish to show there are three distinct entries:  
$ \sigma^{s}{(l-1)}\neq  n-2l$,
$\sigma^{s}{(l-1)}\neq  n-2(m+1)$,
$ n-2l\neq  n-2(m+1)$.
Suppose $ n-2l=n-2(m+1)\Rightarrow l= m+1$, which is a contradiction.
Suppose $ \sigma^{s}{(l-1)}\equiv  n-2l$ and we wish to show there are three distinct entries:  
$n-2(m+1) \neq  n-2l$,
$ \sigma^{s}{(m-1)}\neq n-2l$,
$\sigma^{s}{(m-1)}\neq  n-2(m+1)$.
Since  $\sigma^{s}{(l-1)}\equiv  n-2l$ and $ \sigma^{s}{(l-1)}\neq  \sigma^{s}{(m-1)}\Rightarrow  \sigma^{s}{(m-1)}\neq n-2l$.
Subtracting two equations  $\sigma^{s}{(m-1)}\equiv  n-2(m+1)$  and  $ \sigma^{s}{(l-1)}\equiv  n-2l\Rightarrow m \equiv l-\frac{2}{3}$ $\pmod {n-1}$, which is a contradiction.
Going back to the assumption $ \sigma^{s}{(l-1)}\equiv  n-2m$ we wish to show there are three distinct entries:  
$n-2(m+1)\neq  n-2l$,
$ \sigma^{s}{(m-1)}\neq n-2(m+1)$,
$\sigma^{s}{(m-1)}\neq n-2l$.
Since $\sigma^{s}{(l-1)}\equiv  n-2(m+1)$ and $ \sigma^{s}{(l-1)}\neq  \sigma^{s}{(m-1)}\Rightarrow  \sigma^{s}{(m-1)}\neq n-2(m+1)$.
Subtracting two equations   $\sigma^{s}{(m-1)}\equiv  n-2l$  and $ \sigma^{s}{(l-1)}\equiv  n-2(m+1) \Rightarrow m \equiv l-2$ $\pmod {n-1}$, which is a contradiction.
One of the following two sets has three distinct elements: 
$\{\sigma^{s}{(l-1)},n-2l, n-2(m+1)\}$ or $\{ \sigma^{s}{(m-1)},n-2l, n-2(m+1)\}$.

{\bf 3.}  $G_3$ with  $2 \leq l \leq \frac{n-4}{6}$, $m=\frac{n}{2}-1$, $2\leq i \leq n-1 $, $j=n$

 $$\begin{pmatrix}
                    \sigma^{s}{(l-1)}& \sigma^{s}{(m-1)}\cr
                   n-2l&n-9\end{pmatrix}$$

We wish to show there are three distinct entries:  
$ n-2l \neq n-9$,
$n-2l \neq \sigma^{s}{(m-1)}$,
$\sigma^{s}{(m-1)}\neq  n-9$.
Suppose $ n-2l=n-9 \Rightarrow 2l= 9$, which is a contradiction.
Suppose $ \sigma^{s}{(m-1)}\equiv  n-9 \Rightarrow s \equiv \frac{n}{2}-7$ $\pmod{n-1}$ or $s=\frac{n}{2}-7$.
Suppose $ \sigma^{s}{(m-1)}\equiv n-2l \Rightarrow s-2 \equiv \frac{n}{2}-2l$ $\pmod{n-1}$ and  since $l \leq \frac{n-4}{6}\Rightarrow \frac{n}{2}-2l >0\Rightarrow s-2=\frac{n}{2}-2l$. We wish to show there are three distinct entries:  
$\sigma^{s}{(l-1)} \neq n-2l$,
$ \sigma^{s}{(l-1)}\neq  n-9$,
$n-2l\neq  n-9$.
Since  $ \sigma^{s}{(m-1)}\equiv n-2l$ and  $\sigma^{s}{(l-1)} \neq \sigma^{s}{(m-1)} \Rightarrow \sigma^{s}{(l-1)}\neq  n-2l$. 
Suppose $ \sigma^{s}{(l-1)}\equiv n-9 \Rightarrow s+l\equiv n-8$ and we assumed $s-2\equiv \frac{n}{2}-2l$. Subtracting them  $\Rightarrow \frac{n}{2}+l-10 \equiv 0$. For type 3 matrix we have $n=4+6k \Rightarrow l\equiv 8-3k$ and since $k \geq 4 \Rightarrow l<0$, which is a contradiction.
Assume that $s = n/2-7$ and we wish to show there are three distinct entries:  
$\sigma^{s}{(l-1)} \neq n-2l$,
$ \sigma^{s}{(l-1)}\neq  n-9$,
$n-2l\neq  n-9$,
Suppose $\sigma^{s}{(l-1)}\equiv n-9$ and since $s = n/2-7 \Rightarrow l=n/2-1$, this is a contradiction since $l \leq (n-4)/6$. 
Suppose $ \sigma^{s}{(l-1)}\equiv n-2l$ and  $s = n/2-7 \Rightarrow 3l\equiv n/2+8$. Since $n/2=2+3k \Rightarrow 3l\equiv 10+3k$, right part is not divisible by $3$ this is a contradiction.
One of the following two sets has three distinct elements:  
$\{\sigma^{s}{(m-1)},n-2l, n-9 \}$ or $\{ \sigma^{s}{(l-1)},n-2l, n-9 \}$.

Exceptions: For $n=16$
replacing $n-9$ with $n-11$ and following the argument above we get that one of the following two sets has three distinct elements:  
$\{\sigma^{s}{(m-1)},n-2l, n-11 \}$ or $\{ \sigma^{s}{(l-1)},n-2l, n-11 \}$. 
For $n=22$ 
replacing $n-9$ with $n-5$ and following the same argument we get that one of the following two sets has three distinct elements: 
$\{\sigma^{s}{(m-1)},n-2l, n-5 \}$ or $\{ \sigma^{s}{(l-1)},n-2l, n-5 \}$.

{\bf 4.}  $G_3$ with  $\frac{n+2}{6} \leq l \leq \frac{n}{2}-2$, $m=\frac{n}{2}-1$, $2\leq i \leq n-1 $, $j=n$, and let $s=(n+1)-i$.

 $$\begin{pmatrix}
                    \sigma^{s}{(l-1)}& \sigma^{s}{(m-1)}\cr
                   n-2(l+1)&n-9\end{pmatrix}$$

 We wish to show there are three distinct entries:  
$ n-2l-2 \neq n-9$,
$n-2l-2 \neq \sigma^{s}{(m-1)}$,
$\sigma^{s}{(m-1)}\neq  n-9$.
Suppose $ n-2l-2=n-9 \Rightarrow 2l= 7$, which is a contradiction.
Suppose $ \sigma^{s}{(m-1)}\equiv  n-9 \Rightarrow s \equiv \frac{n}{2}-7$ $\pmod{n-1}$ or $s=\frac{n}{2}-7$.
Suppose $ \sigma^{s}{(m-1)}\equiv n-2l-2 \Rightarrow s\equiv \frac{n}{2}-2l$ $\pmod{n-1}$ and  since $l \leq \frac{n-4}{6}\Rightarrow \frac{n}{2}-2l>0\Rightarrow s=\frac{n}{2}-2l$. We wish to show there are three distinct entries:  
$\sigma^{s}{(l-1)} \neq n-2l-2$,
$ \sigma^{s}{(l-1)}\neq  n-9$,
$n-2l-2\neq  n-9$.
Since  $ \sigma^{s}{(m-1)}\equiv n-2l-2$ and  $\sigma^{s}{(l-1)} \neq \sigma^{s}{(m-1)} \Rightarrow \sigma^{s}{(l-1)}\neq  n-2l-2$. 
Suppose $ \sigma^{s}{(l-1)}\equiv n-9 \Rightarrow s+l\equiv n-8$ and we assumed $s\equiv \frac{n}{2}-2l$. Subtracting them  $\Rightarrow \frac{n}{2}+l-8 \equiv 0$. For type 3 matrix we have $n=4+6k \Rightarrow l\equiv 6-3k$ and since $k \geq 4 \Rightarrow l<0$, which is a contradiction.
Assume that $s = n/2-7$ and we wish to show there are three distinct entries:  
$\sigma^{s}{(l-1)} \neq n-2l-2$,
$ \sigma^{s}{(l-1)}\neq  n-9$,
$n-2l-2\neq  n-9$,
Suppose $\sigma^{s}{(l-1)}\equiv n-9$ and since $s = n/2-7 \Rightarrow l=n/2-1$, this is a contradiction since $l \leq (n-4)/6$. 
Suppose $ \sigma^{s}{(l-1)}\equiv n-2l-2$ and  $s = n/2-7 \Rightarrow 3l\equiv n/2+6$. Since $n/2=2+3k \Rightarrow 3l\equiv 8+3k$, right part is not divisible by $3$ this is a contradiction.
One of the following two sets has three distinct elements:  
$\{\sigma^{s}{(m-1)},n-2l-2, n-9 \}$ or $\{ \sigma^{s}{(l-1)},n-2l-2, n-9 \}$.

Exceptions:
For $n=16$
replacing $n-9$ with $n-11$ and following the argument above we get that one of the following two sets has three distinct elements:  
$\{\sigma^{s}{(m-1)},n-2l-2, n-11 \}$ or $\{ \sigma^{s}{(l-1)},n-2l-2, n-11 \}$. 
For $n=22$ 
replacing $n-9$ with $n-5$ and following the same argument we get that one of the following two sets has three distinct elements: 
$\{\sigma^{s}{(m-1)},n-2l-2, n-5 \}$ or $\{ \sigma^{s}{(l-1)},n-2l-2, n-5 \}$.                 
 
\vspace{.1in}

\emph{\textbf{Case 2}}: We take matrix $G(i,j;l,m)$ with ($2\leq l \leq \frac{n}{2}-1$, $\frac{n}{2} \leq m \leq n-2$), or ($\frac{n}{2} \leq l < m \leq n-2$) and $2\leq i \leq n-1$, $j=n$. Three distinct entries are $\{\sigma^{s}{(l-1)}, \sigma^{s}{(m-1)},n\}$

\vspace{.1 in}

\emph{\textbf{Case 3}}: We take matrix $G(i,j;l,m)$ with $2\leq l\leq \frac{n}{2}-1$, $ m = n-1$, $2\leq i \leq n-1 $, $j=n$.

\noindent \textit{Consider  $G_1$}.

$$\begin{pmatrix}
                  \sigma^{s}{(l-1)}&\sigma^{s}{(m-1)}\cr
                   n-2l&n-1\end{pmatrix}$$

We wish to show there are three distinct entries:
$n-2l \neq n-1$,
$n-2l \neq \sigma^{s}{(n-2)}$,
$\sigma^{s}{(n-2)}\neq n-1$.
Suppose  $n-2l = n-1\Rightarrow l=\frac{1}{2}$, which is a contradiction.
Suppose  $n-2l \equiv \sigma^{s}{(n-2)}\Rightarrow n-2l \equiv s+n-2$ $\pmod {n-1}\Rightarrow s+2l\equiv 2$.
Suppose  $n-1 \equiv \sigma^{s}{(n-2)}\Rightarrow s\equiv 1 \Rightarrow i\equiv n $, which is a contradiction.
Assume $s+2l\equiv 2$ and we wish to show there are three distinct entries:
$\sigma^{s}{(l-1)} \neq \sigma^{s}{(n-2)}$,
$ \sigma^{s}{(l-1)} \neq n-1$,
$\sigma^{s}{(n-2)}\neq n-1$.
Suppose   $ \sigma^{s}{(n-2)}\equiv  n-1$ $\pmod {n-1}\Rightarrow s\equiv 1$, which is a contradiction.
Suppose  $\sigma^{s}{(l-1)} \equiv n-1\Rightarrow s+l\equiv n$ $\pmod {n-1}$.
Subtracting two equations $s+2l\equiv 2$ and $s+l\equiv n \Rightarrow l\equiv 2-n$ $\pmod {n-1}$, which is a contradiction.
One of the following two sets has three distinct elements:
$\{\sigma^{s}{(m-1)},n-2l, n-1\}$ or $\{\sigma^{s}{(l-1)},\sigma^{s}{(m-1)}, n-1\}$.

\vspace{.1in}

\noindent \textit{Consider $G_2$} with $l \leq  \frac{n}{2}-1$, $m =n-1$, $2\leq i \leq n-1 $, $j=n$.

{\bf 1.} $G_2$ with $l\leq \frac{n}{2}-2$. 

$$\begin{pmatrix}
                  \sigma^{s}{(l-1)}&\sigma^{s}{(m-1)}\cr
                   n-2(l+1)&n-1\end{pmatrix}$$

We wish to show there are three distinct entries:
$\sigma^{s}{(l-1)}\neq \sigma^{s}{(n-2)}$,
$\sigma^{s}{(l-1)} \neq  n-2(l+1)$,
$\sigma^{s}{(n-2)}\neq  n-2(l+1)$.
Suppose  $n-2l \equiv \sigma^{s}{(n-2)}\Rightarrow n-2l \equiv s+n-2$ $\pmod {n-1}\Rightarrow s\equiv -2l+2$ $\pmod {n-1}$, which is a contradiction.
Suppose  $\sigma^{s}{(n-2)}\equiv  n-2(l+1)\Rightarrow s\equiv -2l$ $\pmod {n-1}$, which is a contradiction.
Suppose $\sigma^{s}{(l-1)} \equiv  n-2(l+1)\Rightarrow s+3l\equiv n-1$ $\pmod {n-1}$ and we wish to show there are three distinct entries:
$\sigma^{s}{(l-1)} \neq \sigma^{s}{(n-2)}$,
$ \sigma^{s}{(n-2)} \neq n-1$,
$\sigma^{s}{(l-1)}\neq n-1$.
Suppose   $ \sigma^{s}{(n-2)}\equiv  n-1$ $\pmod {n-1}$ $\Rightarrow s \equiv 1 \Rightarrow i\equiv n $, which is a contradiction.
Suppose $ \sigma^{s}{(l-1)} \equiv  n-1$ $\pmod {n-1}$ $\Rightarrow s+l \equiv n$.
Subtracting two equations $s+3l\equiv n-1$ and $s+l\equiv n \Rightarrow 2l\equiv -1$ $\pmod {n-1}$, which is a contradiction.
One of the following two sets has three distinct elements: 
$\{\sigma^{s}{(l-1)},\sigma^{s}{(m-1)},n-2(l+1)\}$ or $\{ \sigma^{s}{(l-1)},\sigma^{s}{(m-1)}, n-1\}$.

{\bf 2.} $G_2$ with $l =\frac{n}{2}-1$.

$$\begin{pmatrix}
                \sigma^{s}{(l-1)}&\sigma^{s}{(m-1)} \cr
                   Y&n-1\end{pmatrix}$$

If  \textit{$K$ is Even} $ \Rightarrow Y=\frac{n}{2}-2$.
We wish to show there are three distinct entries:
$\sigma^{s}{(\frac{n}{2}-2)}\neq \sigma^{s}{(n-2)}$,
$\sigma^{s}{(\frac{n}{2}-2)} \neq  \frac{n}{2}-2$,
$\sigma^{s}{(n-2)}\neq  \frac{n}{2}-2$.
First two inequalities are correct.
$\sigma^{s}{(n-2)}\equiv \frac{n}{2}-2\Rightarrow s+n-2\equiv \frac{n}{2}-2$ $\pmod {n-1}\Rightarrow s\equiv -\frac{n}{2}$, which is a contradiction.
There are three distinct elements $\{\sigma^{s}{(l-1)},\sigma^{s}{(m-1)},Y\}$

\vspace{.1in}

If  \textit{$K$ is Odd} $\Rightarrow  Y=\frac{n}{2}+1$.
We wish to show there are three distinct entries:
$\sigma^{s}{(n-2)}\neq n-1$.
$\sigma^{s}{(n-2)} \neq  \frac{n}{2}+1$,
$\frac{n}{2}+1 \neq  n-1$.
Suppose  $\sigma^{s}{(n-2)}\equiv n-1 \Rightarrow s+n-2\equiv n-1\Rightarrow s \equiv 1$ $\pmod {n-1} \Rightarrow i\equiv n $, which is a contradiction.
Suppose  $\frac{n}{2}+1= n-1\Rightarrow n= 4$, which is a contradiction.
Suppose  $\sigma^{s}{(n-2)}\equiv \frac{n}{2}+1\Rightarrow s+\frac{n}{2}\equiv 3$ $\pmod {n-1}$. Since $s=(n+1)-i \Rightarrow  i\equiv 1.5n-4$, this is a contradiction.
There are three distinct elements $\{\sigma^{s}{(m-1)},n-1,Y\}$.

\vspace{.1in}

\noindent \textit{Consider $G_3$} with $l \leq  \frac{n}{2}-1$, $m =n-1$, $2\leq i \leq n-1 $, $j=n$.

{\bf 1.} $G_3$  with $2\leq l \leq \frac{n-4}{6}$.

$$\begin{pmatrix}
                  \sigma^{s}{(l-1)}& \sigma^{s}{(m-1)}\cr
                   n-2l&n-1\end{pmatrix}$$

We wish to show there are three distinct entries:
$\sigma^{s}{(l-1)}\neq \sigma^{s}{(n-2)}$,
$\sigma^{s}{(l-1)} \neq  n-1$,
$\sigma^{s}{(n-2)}\neq  n-1$.
Suppose  $\sigma^{s}{(n-2)}\equiv  n-1\Rightarrow s+n-2\equiv n-1\Rightarrow i\equiv n$ $\pmod {n-1}$, which is a contradiction.
Suppose  $\sigma^{s}{(l-1)}\equiv n-1\Rightarrow s+l\equiv n$ $\pmod {n-1}$ and we wish to show there are three distinct entries:
$\sigma^{s}{(n-2)} \neq n-1$,
$n-2l \neq n-1$,
$\sigma^{s}{(n-2)}\neq n-2l$.
Suppose   $n-2l= n-1\Rightarrow l=\frac{1}{2}$, which is a contradiction.
Suppose  $\sigma^{s}{(n-2)} \equiv n-2l\Rightarrow s+2l=2$ $\pmod {n-1}$.
Subtracting two equations $s+2l=2$ and $s+l\equiv n\Rightarrow  l\equiv 2-n$ $\pmod {n-1}$, which is a contradiction.
One of the following two sets has three distinct elements:  
$\{\sigma^{s}{(l-1)},\sigma^{s}{(m-1)},n-1\}$ or  $\{\sigma^{s}{(m-1)},n-2l,n-1\}$.

{\bf 2.} $G_3$  with $ \frac{n-4}{6} \leq l \leq \frac{n}{2}-2 $.

$$\begin{pmatrix}
                   \sigma^{s}{(l-1)}& \sigma^{s}{(m-1)}\cr
                   n-2(l+1)&n-1\end{pmatrix}$$

We wish to show there are three distinct entries:
$\sigma^{s}{(l-1)}\neq \sigma^{s}{(n-2)}$,
$\sigma^{s}{(n-2)} \neq  n-1$,
$\sigma^{s}{(l-1)}\neq  n-1$.
Suppose  $\sigma^{s}{(n-2)}\equiv  n-1\Rightarrow i \equiv n$ $\pmod {n-1}$, which is a contradiction.
Suppose  $\sigma^{s}{(l-1)}\equiv n-1\Rightarrow s+l\equiv n$ $\pmod {n-1}$ and we wish to show there are three distinct entries:
$\sigma^{s}{(l-1)} \neq \sigma^{s}{(n-2)}$,
$\sigma^{s}{(l-1)} \neq n-2(l+1)$,
$\sigma^{s}{(n-2)}\neq n-2(l+1)$.
Suppose   $ \sigma^{s}{(n-2)} \equiv n-2(l+1)$ $\pmod {n-1}\Rightarrow s\equiv -2l$  $\pmod {n-1}$, which is a contradiction.
Suppose   $\sigma^{s}{(l-1)} \equiv n-2(l+1)\Rightarrow s+3l \equiv n-1$  $\pmod {n-1}$.
Subtracting this equation from $s+l\equiv n$ $\pmod {n-1}\Rightarrow 2l\equiv -1$, which is a contradiction.
One of the following two sets has three distinct elements:
 $\{\sigma^{s}{(l-1)},\sigma^{s}{(m-1)},n-1\}$ or $\{\sigma^{s}{(m-1)},\sigma^{s}{(l-1)},n-2(l+1)\}$.

{\bf 3.} $G_3$  with   $l=\frac{n}{2}-1$.

 $$\begin{pmatrix}
                  \sigma^{s}{(l-1)}& \sigma^{s}{(m-1)}\cr
                   n-9&n-1\end{pmatrix}$$

We wish to show there are three distinct entries:
$n-9\neq n-1$,
$\sigma^{s}{(n-2)} \neq n-1$,
$\sigma^{s}{(n-2)}\neq  n-9$.
Suppose  $\sigma^{s}{(n-2)}\equiv  n-1\Rightarrow s \equiv 1$ $\pmod {n-1}$, which is a contradiction.
Suppose  $\sigma^{s}{(n-2)}\equiv n-9\Rightarrow s\equiv -7$ $\pmod {n-1}$, which is a contradiction.
There are three distinct elements $\{\sigma^{s}{(m-1)},n-1, n-9\}$.

Exceptions:
For $n=16$ three distinct elements are: $\{\sigma^{s}{(7)},5,15\}$.
For $n=22$ three distinct elements are: $\{\sigma^{s}{(10)},21,17\}$.

\vspace{.1 in}

\emph{\textbf{Case 4}}: We take matrix $G(i,j;l,m)$ with  $\frac{n}{2} \leq l \leq n-2$, ($m=n-1$ or $m=n$), $2\leq i \leq n-1$, $j=n$. Three distinct entries are $\{\sigma^{s}{(l-1)},\sigma^{s}{(m-1)}, n\}$.

\vspace{.1 in}

\emph{\textbf{Case 5}}:  We take matrix $G(i,j;l,m)$ with $l=n-1$, $m=n$, $2\leq i \leq n-1 $, $j=n$. Three distinct entries are  $\{n-1,\sigma^{s}{(l-1)},\sigma^{s}{(m-1)}\}$

\vspace{.2in}

\noindent \emph{\textbf{Step 4}}:

For matrix $G(i,j;l,m)$ with $i=1$, $2 \leq j \leq \n$,  $l=1$, $2\leq m \leq n$ we look at four cases and consider every matrix type defined above of even $n\pmod 6$.

\vspace{.1 in}

\emph{\textbf{Case 1}}: We take matrix $G(i,j;l,m)$ with $l=1$, $2\leq m\leq n$, $i =1$, $j = 2$.

$$\begin{pmatrix}
                    1& m\cr
                    3& m-1\end{pmatrix}$$

We know that $m \neq 1$ and $m \neq m-1$. If $m-1 \neq 1\Rightarrow $  we have three distinct entries.
Suppose $m-1=1\Rightarrow m=2 \Rightarrow$ we have three distinct entries.
One of the following two sets has three distinct elements:  $\{m,1,m-1\}$ or $\{1,3,m\}$.

\vspace{.1 in}

\emph{\textbf{Case 2}}: We take matrix $G(i,j;l,m)$ with $l =1$, $2\leq m\leq n$, $i =1$, $3 \leq j \leq \frac{n}{2}+1$.

$$\begin{pmatrix}
                    1& m\cr
                    n& X\end{pmatrix}$$

There are three distinct entries $\{1, m, n\}$.

\vspace{.1 in}

\emph{\textbf{Case 3}}: We take matrix  $G(i,j;l,m)$ with $l=1$, $2\leq m\leq n$, $i =1$, $\frac{n}{2}+2 \leq j \leq n-1$, and let $t=(n+1)-j\Rightarrow t=\frac{n}{2}-1$.
\vspace{.1in}

\noindent \textit{Consider  $G_1$} and let $t=(n+1)-j$.

$$\begin{pmatrix}
                    1& m\cr
                    2(j-1)-n& \sigma^{t}{(m-1)}\end{pmatrix}$$

We wish to show there are three distinct entries:
$m \neq 1$,
$m \neq \sigma^{t}{(m-1})$,
$1 \neq \sigma^{t}{(m-1)}$.
Suppose  $1 \equiv \sigma^{(t)}{(m-1)}\Rightarrow t \equiv 2-n$. Hence $t=(n+1)-j \Rightarrow n-1 \equiv j-m$ $\pmod {n-1}$ and we wish to show there are three distinct entries:
$m \neq 1$,
$2j-2-n \neq 1$,
$m\neq 2j-2-n$.
Suppose $2j-2-n=1\Rightarrow j=\frac{n}{2}-\frac{1}{2}$. Since $j >\frac{n}{2}+1$, this is a contradiction.
Suppose $m =2j-2-n$,
subtracting it from  $n-1\equiv j-m  \Rightarrow j \equiv 3$, which is a contradiction.
One of the following two sets has three distinct elements:  $\{m,\sigma^{t}{(m-1)},1\}$ or $\{m,2j-2-n,1\}$ 
\vspace{.1in}

\noindent \textit{Consider  $G_2$}.

{\bf1.}  $G_2$ with    $j=\frac{n}{2}+2$.

$$\begin{pmatrix}
                    1& m\cr
                  Y& \sigma^{t}{(m-1)}\end{pmatrix}$$

If  \textit{$K$ is Even} $\Rightarrow Y=\frac{n}{2}-2$.
We wish to show there are  three distinct entries: 
$m \neq 1$,
$\frac{n}{2}-2 \neq 1$,
$\frac{n}{2}-2 \neq m$.
Suppose  $\frac{n}{2}-2 = 1\Rightarrow n=6$, which is a contradiction.
Suppose  $\frac{n}{2}-2 = m \Rightarrow  n=2m+4$ and we wish to show there are  three distinct entries: 
$\sigma^{\frac{n}{2}-1}{(m-1)} \neq m$,
$m \neq 1$,
$\sigma^{\frac{n}{2}-1}{(m-1)}\neq 1$.
Suppose $\sigma^{\frac{n}{2}-1}{(m-1)}\equiv 1\Rightarrow  n \equiv 6-2m$.
Subtracting it from e  $n \equiv 6-2m \Rightarrow  m \equiv \frac{1}{2}$ $\pmod {n-1}$, which is a contradiction.
One of the following two sets has three distinct elements:  
$\{m,Y,1\}$ or $\{m,\sigma^{t}{(m-1)},1\}$ 

\vspace{.1in}

If  \textit{$K$ is Odd} $ \Rightarrow Y=\frac{n}{2}+1$.
We wish to show there are  three distinct entries:
$m \neq 1$,
$\frac{n}{2}+1\neq 1$,
$\frac{n}{2}+1 \neq m$.
First two inequalities are correct.
Suppose  $\frac{n}{2}+1 = m \Rightarrow  n=2m-2$ and we wish to show there are  three distinct entries:
$\sigma^{\frac{n}{2}-1}{(m-1)} \neq m$,
$m \neq 1$,
$\sigma^{\frac{n}{2}-1}{(m-1)}\neq 1$.
Suppose $\sigma^{\frac{n}{2}-1}{(m-1)}\equiv 1\Rightarrow  n \equiv 2-2m$ $\pmod {n-1}$.
Subtracting it from $ n=2m-2 \Rightarrow m \equiv 1$, which is a contradiction.
One of the following two sets has three distinct elements:    
 $\{m,Y,1\}$ or $\{m,\sigma^{t}{(m-1)},1\}$ 

{\bf 2.} $G_2$ with $\frac{n}{2}+3 \leq j \leq n-1$.

$$\begin{pmatrix}
                    1& m\cr
                    2(j-2)-n& \sigma^{t}{(m-1)}\end{pmatrix}$$

We wish to show there are three distinct entries:
$m \neq 1$,
$m\neq \sigma^{t}{(m-1)}$,
$\sigma^{t}{(m-1)}\neq 1$.
Suppose $\sigma^{t}{(m-1)}\equiv 1\Rightarrow j \equiv n-1+m$ $\pmod {n-1}$. Hence $j \leq n-1$, this is a contradiction
There are three distinct entries $\{m,\sigma^{t}{(m-1)},1)\}$. 
\vspace{.1in}

\noindent \textit{Consider  $G_3$}.

{\bf 1.} $G_3$ with  $j=\frac{n}{2}+2$.

$$\begin{pmatrix}
                    1& m\cr
                    n-9& \sigma^{t}{(m-1)}\end{pmatrix}$$

We wish to show there are  three distinct entries:
$m \neq 1$,
$\sigma^{t}{(m-1)}\neq m$,
$\sigma^{t}{(m-1)}\neq 1$.
First two inequalities are correct.
Suppose  $\sigma^{t}{(m-1)}\equiv 1\Rightarrow m \equiv 2-\frac{n}{2}$ $\pmod {n-1}$, which is a contradiction.
Three distinct elements:  $\{m,\sigma^{t}{(m-1)},1\}$.

Exceptions:
For $n=16$ following the argument above we have that three distinct elements
$\{m,\sigma^{t}{(m-1)},1\}$.
For $n=22$ three distinct entries are $\{m,\sigma^{t}{(m-1)},1\}$.

{\bf 2.} $G_3$ with $\frac{n}{2}+3 \leq j \leq \frac{5n+4}{6}$. 

$$\begin{pmatrix}
                    1& m\cr
                    2(j-2)-n& \sigma^{t}{(m-1)}\end{pmatrix}$$

We wish to show there are  three distinct entries:
$m \neq 1$,
$\sigma^{t}{(m-1)}\neq m$,
$\sigma^{t}{(m-1)}\neq 1$.
First two inequalities are correct. Suppose  $\sigma^{t}{(m-1)}\equiv 1\Rightarrow n-j+m \equiv 1$ $\pmod {n-1}$ and we wish to show there are  three distinct entries:
$m \neq 1$,
$2(j-2)-n \neq 1$,
$2(j-2)-n\neq m$.
Suppose $ 2(j-2)-n=1 \Rightarrow j=\frac{n+5}{2} \Rightarrow j=\frac{n}{2}+2.5$ and since $\frac{n}{2}+3 \leq j$, this is a contradiction.
Suppose  $2(j-2)-n=m \Rightarrow 2j=n+m+4$. Subtracting this equation from  $ n-j+m \equiv 1 \Rightarrow j \equiv 5$ $\pmod {n-1}$. Hence  $n\geq 16$ and $j\leq \frac{5n+4}{6}\Rightarrow j\geq 14$, which is a contradiction. 
One of the following two sets has three distinct elements:  
$\{m,\sigma^{t}{(m-1)},1\}$ or $\{m, 1, 2(j-2)-n\}$.

{\bf 3.} $G_3$ with  $\frac{5n+10}{6}\leq j \leq n-1$.

$$\begin{pmatrix}
                    1& m\cr
                    2(j-1)-n& \sigma^{t}{(m-1)}\end{pmatrix}$$

We wish to show there are three distinct entries:
$m \neq 1$,
$\sigma^{t}{(m-1)}\neq m$,
$\sigma^{t}{(m-1)}\neq 1$.
First two inequalities are correct. 
Suppose  $\sigma^{t}{(m-1)}\equiv 1\Rightarrow n-j+m \equiv 1$ $\pmod {n-1}$ and we wish to show there are three distinct entries:
$m \neq 1$,
$2(j-1)-n \neq 1$,
$2(j-1)-n \neq m$.
Suppose $2(j-1)-n=1 \Rightarrow j=\frac{n}{2}+1.5$ and since $\frac{5n+10}{6}\leq j$, this is a contradiction.
Suppose  $2(j-1)-n=m$, subtracting it from  $m \equiv j+1-n \Rightarrow j \equiv 3$. Hence $n\geq 16$ and $j \geq \frac{5n+10}{6}\Rightarrow  j\geq 15$, which is a contradiction.
One of the following two sets has three distinct elements:  
 $\{m,\sigma^{t}{(m-1)},1\}$ or $\{m, 1, 2(j-1)-n\}$.

\vspace{.1 in}

\emph{\textbf{Case 4}}:  We take matrix $G(i,j;l,m)$ with $l=1$, $2 \leq m \leq n$, $i=1$, $j=n$,

\noindent \textit{Consider  $G_1$} with $2\leq m \leq  \frac{n}{2}-1$ and \textit {$G_3$} with $2\leq m \leq \frac{n-4}{6}$

$$\begin{pmatrix}
                    1& m\cr
                    n-2&n-2m\end{pmatrix}$$

We wish to show there are three distinct entries:
$n-2\neq 1$,
$m \neq  1$ ,
$m\neq  n-2$.
Suppose  $m=n-2\Rightarrow n=m+2$ and wish to show there are three distinct entries:
$n-2\neq n-2m$,
$n-2 \neq  1$,
$n-2m\neq 1$.
First two inequalities are correct.
Suppose  $n-2m=1$, subtracting it from $n=m+2\Rightarrow m=1$, which is a contradiction.
One of the following two sets has three distinct elements:  $\{m,n-2,1\}$ or $\{n-2,1,n-2m\}$.

\vspace{.1 in}
 
\noindent \textit{Consider $G_2$} with  $2\leq m \leq  \frac{n}{2}-1$ and \textit{$G_3$} with $\frac{n+2}{6}\leq m \leq \frac{n}{2}-1$.

{\bf 1.} $G_2$ with  $2\leq m \leq  \frac{n}{2}-2$ and  $G_3$ with $\frac{n+2}{6}\leq m \leq \frac{n}{2}-2$.

$$\begin{pmatrix}
                    1& m\cr
                    n-2&n-2(m+1)\end{pmatrix}$$

We wish to show there are three distinct entries:
$n-2\neq 1$,
$m \neq  1$,
$m \neq  n-2$.
Assume $m=n-2$ and  wish to show there are three distinct entries:
$n-2\neq n-2(m+1)$,
$n-2 \neq  1$,
$n-2(m+1)\neq 1$.
Suppose  $n-2(m+1)=1$, subtracting it from $n=m+2$  $\Rightarrow m=-1$, which is a contradiction.
One of the following two sets has three distinct elements:  $\{m,n-2,1\}$ or $\{n-2,1,n-2(m+1)\}$.

{\bf 2.} $G_2$ and $G_3$ with $m=\frac{n}{2}-1$.

$$\begin{pmatrix}
                    1& \frac{n}{2}-1\cr
                    n-2&Y\end{pmatrix}$$
Three distinct elements are $\{1, n-2, \frac{n}{2}-1\}$.
\vspace{.1in}

\noindent \textit{Consider  $G_1$, $G_2$, $G_3$}

{\bf 1.} $G(i,j;l,m)$ with  ($ \frac{n}{2} \leq m \leq n-2$ or $m=n$).
Three distinct entries are $\{1,n-2,n)\}$.

{\bf 2.} $G(i,j;l,m)$ with $m=n-1$.
Three distinct entries are $\{1,n-1,n-2\}$.

\vspace{.1in}

\noindent \emph{\textbf{Step 5}}:

For matrix $G(i,j;l,m)$ with $i=2$, $3 \leq j \leq \n$,  $l=1$, $2\leq m \leq n$ we look at three cases and consider every matrix type defined above of even $n\pmod 6$.
\vspace{.1in}

\emph{\textbf{Case 1}}: We take the submatrix of $G(i,j;l,m)$ with $l =1$, $2\leq m\leq n$, $i =2$, $3 \leq j \leq \frac{n}{2}+1$.

$$\begin{pmatrix}
                    3& m-1\cr
                    n&\sigma^{t}{(m-1)}\end{pmatrix}$$

\noindent There are three distinct entries $\{m-1, \sigma^{t}{(m-1}), n\}$.

\vspace{.1 in}

\emph{\textbf{Case 2}}:  We take the submatrix of $G(i,j;l,m)$ with $l =1$, $2\leq m\leq n$, $i =2$, $\frac{n}{2}+2 \leq j \leq n-1$, and let  $t=(n+1)-j$
\vspace{.1in}

\noindent \textit{Consider  $G_1$} 

$$\begin{pmatrix}
                    3& m-1\cr
                    2(j-1)-n& \sigma^{t}{(m-1)}\end{pmatrix}$$

\noindent We wish to show there are  three distinct entries:
$m -1  \neq \sigma^{t}{(m-1})$,
$2j-n-2 \neq \sigma^{t}{(m-1})$,
$2j-n-2  \neq m-1$.
The first inequality  is correct.
Suppose  $2j-n-2  \equiv \sigma^{(t)}{(m-1)}\Rightarrow  3j \equiv 2n+2+m$ $\pmod {n-1}$.
Suppose  $2j-n-2=m-1\Rightarrow 2j= n+1+m$ and we wish to show there are  three distinct entries:
$m-1 \neq 3$,
$m-1 \neq \sigma^{t}{(m-1})$,
$\sigma^{t}{(m-1}) \neq 3$.
Suppose $m-1 = 3\Rightarrow m=4$.
Suppose $\sigma^{t}{(m-1}) \equiv 3\Rightarrow j  \equiv n+m-3$ $\pmod {n-1}$.
Subtracting two equations $m=4$ and $3j \equiv 2n+2+m\Rightarrow 3j \equiv 2n+6$ and hence $n=2+6k \Rightarrow 3j\equiv 2(6k+5)$, right part is not divisible by 3, which is a contradiction.
Subtracting two equations $3j= 2n+2+m$ and $j \equiv n+m-3\Rightarrow  2j \equiv n+5$, and since $n$ is even, $j$ is not integer,  which is a contradiction.
Subtracting two equations $m=4$ and $2j=n+1+m \Rightarrow  2j =n+5$, which is a contradiction.
Subtracting two equations $2j= n+1+m$ and $j \equiv n+m-3\Rightarrow j \equiv 4$ and since  $j \geq \frac{n}{2}+2$, this is a contradiction.
One of the following two sets has three distinct elements:  
$\{m-1,\sigma^{t}{(m-1}),2j-n-2\}$ or $\{m-1,3,\sigma^{t}{(m-1})\}$.
\vspace{.1in}

\noindent \textit{Consider $G_2$}

{\bf 1.}  $G_2$ with $l =1$, $2\leq m\leq n$, $i =2$, $j=\frac{n}{2}+2$, $t=(n+1)-j\Rightarrow t=\frac{n}{2}-1$.

$$\begin{pmatrix}
                    3& m-1\cr
                    Y& \sigma^{t}{(m-1)}\end{pmatrix}$$

If  \textit{$K$ is Even} $ \Rightarrow Y=\frac{n}{2}-2$.
We wish to show there are three distinct entries:
$m-1  \neq \sigma^{t}{(m-1})$,
$\frac{n}{2}-2 \neq \sigma^{t}{(m-1})$,
$\frac{n}{2}-2  \neq m-1$.
The first inequality is correct.
Suppose  $\frac{n}{2}-2 \equiv \sigma^{t}{(m-1})\Rightarrow \frac{n}{2}-2  \equiv \frac{n}{2}-1+m-1\Rightarrow m \equiv 0$ $\pmod {n-1}$, which is a contradiction.
Suppose  $\frac{n}{2}-2 = m-1\Rightarrow \frac{n}{2}-1=m$ and we wish to show there are three distinct entries:
$m-1 \neq 3$,
$m-1 \neq \sigma^{t}{(m-1})$,
$\sigma^{t}{(m-1}) \neq 3$.
Suppose $\sigma^{t}{(m-1}) \equiv 3\Rightarrow \frac{n}{2}-1+m-1  \equiv 3\Rightarrow n-10 \equiv -2m$ $\pmod {n-1} \Rightarrow  n < 10$, which is a contradiction. 
Suppose $m-1=3\Rightarrow m=4$ and since  $\frac{n}{2}-1=m\Rightarrow n=10$, this is  a contradiction.
One of the following two sets has three distinct elements:  
$\{m-1,\sigma^{t}{(m-1}), Y\}$ or $\{m-1, 3, \sigma^{t}{(m-1})\}$.

\vspace{.1in}

If  \textit{$K$ is Odd} $ \Rightarrow Y=\frac{n}{2}+1$, and let  $t=(n+1)-j \Rightarrow t=\frac{n}{2}-1$.
We wish to show there are three distinct entries:
$m-1  \neq \sigma^{t}{(m-1})$,
$\frac{n}{2}+1 \neq \sigma^{t}{(m-1})$,
$\frac{n}{2}+1  \neq m-1$.
The first inequality is correct.
Suppose  $\frac{n}{2}+1 \equiv \sigma^{t}{(m-1})\Rightarrow  m \equiv 3$ $\pmod {n-1}$.
Suppose  $\frac{n}{2}+1 = m-1\Rightarrow n=2m-4$ and we wish to show there are three distinct entries:
$m-1 \neq 3$,
$m-1 \neq \sigma^{t}{(m-1})$,
$\sigma^{t}{(m-1}) \neq 3$.
From the argument above (for even $k$) we have $\sigma^{t}{(m-1}) \neq 3$ is correct.
Suppose $m-1=3\Rightarrow m=4$ and since $m\equiv 3$ we have a contradiction.
Taking two equations $n=2m-4$ and $m=4\Rightarrow n=4$ we get a contradiction.
One of the following two sets has three distinct elements:
$\{m-1,\sigma^{t}{(m-1}),Y\}$ or $\{m-1,3,\sigma^{t}{(m-1})\}$.

{\bf 2.}  $G_2$ with $l =1$, $2\leq m\leq n$, $i =2$, $\frac{n}{2}+3 \leq j \leq n-1$, and let $t=(n+1)-j$.

$$\begin{pmatrix}
                    3& m-1\cr
                    2(j-2)-n& \sigma^{t}{(m-1)}\end{pmatrix}$$

We wish to show there are three distinct entries:
$m-1  \neq \sigma^{t}{(m-1})$,
$2(j-2)-n \neq \sigma^{t}{(m-1})$,
$2(j-2)-n  \neq m-1$.
The first inequality is correct.
Suppose  $2(j-2)-n \equiv \sigma^{t}{(m-1})\Rightarrow 3j-4-2n \equiv m$  $\pmod {n-1}$.
Suppose  $2(j-2)-n = m-1\Rightarrow 2j-n-3=m$ and we wish to show there are three distinct entries:
$m-1 \neq 3$,
$m-1 \neq \sigma^{t}{(m-1})$,
$\sigma^{t}{(m-1}) \neq 3$.
Suppose $m-1=3\Rightarrow m=4$ and since $3j-4-2n \equiv m \Rightarrow 3j\equiv  2n+8$ and since $n=6+6k \Rightarrow j\equiv  4k+6\frac{2}{3}$ which is a contradiction.
Now take  $2j-n-3=m$ and $m=4\Rightarrow j=\frac{n}{2}+2.5$, which is a contradiction.
Suppose $\sigma^{t}{(m-1}) \equiv 3\Rightarrow m  \equiv 3+j-n$  $\pmod {n-1}$ and since  $3j-4-2n \equiv m \Rightarrow 2j\equiv n+7$  $\pmod {n-1}\Rightarrow j\equiv \frac{n}{2}+3.5$, which is a contradiction.
Now take   $2j-n-3=m \Rightarrow j=6$. Since $j \geq \frac{n}{2}+3$, this is a contradiction.
One of the following two sets has three distinct elements:  
$\{m-1,\sigma^{t}{(m-1}),2(j-1)-n\}$ or $\{m-1,3,\sigma^{t}{(m-1})\}$.

\vspace{.1in}

\noindent \textit{Consider $G_3$}

{\bf 1.} $G_3$ with $l =1$, $2\leq m\leq n$, $i=2$, $j=\frac{n}{2}+2$, and let $t=(n+1)-j\Rightarrow t=\frac{n}{2}-1$.

$$\begin{pmatrix}
                    3& m-1\cr
                    n-9& \sigma^{t}{(m-1)}\end{pmatrix}$$

We wish to show there are three distinct entries:
$m -1  \neq \sigma^{t}{(m-1})$,
$n-9 \neq \sigma^{t}{(m-1})$,
$n-9 \neq m-1$.
The first inequality is correct.
Suppose  $n-9 \equiv \sigma^{t}{(m-1})\Rightarrow n\equiv 2m+14$.
Suppose  $n-9 = m-1\Rightarrow n=m+8$.
We wish to show there are three distinct entries:
$m-1 \neq \sigma^{t}{(m-1})$,
$m-1 \neq 3$,
$\sigma^{t}{(m-1}) \neq 3$.
Suppose $m-1=3\Rightarrow m=4$ and since $n\equiv 2m+14\Rightarrow n\equiv 22$, which is a contradiction.
Now take $n=m+8$ and since $m=4\Rightarrow n=12$, which is a contradiction.
Suppose $\sigma^{t}{(m-1}) \equiv 3\Rightarrow n+2m \equiv 10$  $\pmod {n-1}$ and since $n\equiv 2m+14\Rightarrow m\equiv -1$, which is a contradiction.
Now take  $n=m+8$   and  since $n+2m \equiv 10\Rightarrow m=\frac{2}{3}$, which is a contradiction.
One of the following two sets has three distinct elements:  
$\{m-1,\sigma^{t}{(m-1}),n-9\}$ or  $\{m-1,3,\sigma^{t}{(m-1})\}$.

Exceptions:
For $n=16$ following the argument above we have three distinct elements: $\{m-1,\sigma^{t}{(m-1}),5\}$ or  $\{m-1,3,\sigma^{t}{(m-1})\}$.
For$n=22$

$$\begin{pmatrix}
                    3& m-1\cr
                    17& \sigma^{10}{(m-1)}\end{pmatrix}$$

We wish to show there are three distinct entries:
$17 \neq 3$,
$\sigma^{10}{(m-1)}\neq 3$,
$\sigma^{10}{(m-1)}\neq 17$.
Suppose  $\sigma^{10}{(m-1)}\equiv 3\Rightarrow m \equiv -6$, which is a contradiction.
Suppose  $\sigma^{10}{(m-1)}\equiv 17\Rightarrow m \equiv 8$.
We wish to show there are three distinct entries:
$17 \neq 3$,
$m-1\neq 3$,
$m-1\neq 17$.
Suppose  $m-1= 3$. Since $m \equiv 8$, this is a contradiction.
Suppose  $m-1=17$. Since $m \equiv 8$, this is a contradiction.
One of the following two sets has three distinct elements:    
$\{17, \sigma^{t}{(m-1)},3\}$ or $\{3, 17, m-1\}$.

{\bf 2.} $G_3$ with $l =1$, $2\leq m\leq n$, $i=2$, $\frac{n}{2}+3 \leq j \leq \frac{5n+4}{6}$, and let $t=(n+1)-j$

$$\begin{pmatrix}
                    3& m-1\cr
                    2(j-2)-n& \sigma^{t}{(m-1)}\end{pmatrix}$$

We wish to show there are three distinct entries:
$m -1  \neq \sigma^{t}{(m-1})$,
$2j-4-n\neq \sigma^{t}{(m-1})$,
$2j-4-n\neq m-1$.
Suppose  $2j-4-n \equiv \sigma^{t}{(m-1})\Rightarrow 2n+m+4 \equiv 3j$  $\pmod {n-1}$.
Suppose  $2j-4-n = m-1\Rightarrow n+m+3=2j$.
We wish to show there are three distinct entries:
$m-1 \neq 3$,
$m-1 \neq \sigma^{t}{(m-1})$,
$\sigma^{t}{(m-1}) \neq 3$.
Suppose $m-1=3\Rightarrow m=4$.
Hence  $2n+m+4 \equiv 3j$ and  $m=4\Rightarrow 2n+8 \equiv 3j$ $\pmod {n-1}$. Since $n=4+6k\Rightarrow j\equiv 4k+5\frac{1}{3}$, which is a contradiction.
Hence $n+m+3=2j$ and  $m=4\Rightarrow  j\equiv \frac{n}{2}+3.5$, which is a contradiction.
Suppose $\sigma^{t}{(m-1}) \equiv 3\Rightarrow j \equiv n+m-3$  $\pmod {n-1}$.
Hence  $2n+m+4 \equiv 3j$  and $j \equiv n+m-3$, subtracting  them $\Rightarrow j\equiv \frac{n}{2}+3.5$, which is a contradiction. 
Hence  $n+m+3=2j$  and $j \equiv n+m-3$, subtracting  them $\Rightarrow j \equiv 6$ $\pmod {n-1}$. Since $j\geq \frac{n}{2}+3$, which is a contradiction.
One of the following two sets has three distinct elements:  
$\{m-1,\sigma^{t}{(m-1}),2j-4-n\}$ or $\{m-1,3,\sigma^{t}{(m-1})\}$.

{\bf 3.}  $G_3$ with $l =1$, $2\leq m\leq n$, $i=2$, $ \frac{5n+10}{6} \leq j \leq n-1$, and let $t=(n+1)-j$.

$$\begin{pmatrix}
                    3& m-1\cr
                    2(j-1)-n& \sigma^{t}{(m-1)}\end{pmatrix}$$

We wish to show there are three distinct entries:
$m -1  \neq \sigma^{t}{(m-1})$,
$2j-2-n\neq \sigma^{t}{(m-1})$,
$2j-2-n\neq m-1$.
Suppose  $2j-2-n \equiv \sigma^{t}{(m-1})\Rightarrow 2n+m+2 \equiv 3j$ $\pmod {n-1}$.
Suppose  $2j-2-n = m-1\Rightarrow n+m+1=2j$.
We wish to show there are three distinct entries:
$m-1 \neq 3$,
$m-1 \neq \sigma^{t}{(m-1})$,
$\sigma^{t}{(m-1}) \neq 3$.
Suppose $m-1 =3\Rightarrow m=4$.
Subtracting  two equations   $2n+m+2 \equiv 3j$ and  $m=4 \Rightarrow 2n+6 \equiv 3j$ $\pmod {n-1}$. Since $j \geq \frac{5n+10}{6}$, this is a contradiction. 
Subtracting  two equations  $n+m+1=2j$ and  $m=4 \Rightarrow j=\frac{n}{2}+2.5$, which is a contradiction.
Suppose $\sigma^{t}{(m-1}) \equiv 3\Rightarrow j \equiv n+m-3$  $\pmod {n-1}$.
Subtracting  two equations  $2n+m+2\equiv 3j$  and $j \equiv n+m-3\Rightarrow j=\frac{n}{2}+2.5$, which is a contradiction.
Subtracting  two equations $n+m+1=2j$  and $j \equiv n+m-3 \Rightarrow j \equiv 4$. Since $j \geq \frac{5n+10}{6}$ and $n\geq 28$, this is a contradiction. 
One of the following two sets has three distinct elements:  
$\{m-1,\sigma^{t}{(m-1}),2j-2-n\}$ or $\{m-1,3,\sigma^{t}{(m-1})\}$.

\vspace{.1 in}

\emph{\textbf{Case 3}}: We take matrix $G(i,j;l,m)$ with $l=1$, $2\leq m\leq n$, $i =2$, $j=n$.

\vspace{.1in}

\noindent \textit{Consider  $G_1$} with $2\leq m \leq  \frac{n}{2}-1$ and \textit {consider $G_3$} with $2\leq m \leq \frac{n-4}{6}$.

$$\begin{pmatrix}
                    3& m-1\cr
                    n-2&n-2m\end{pmatrix}$$

We wish to show there are three distinct entries:
$n-2\neq 3$,
$m-1\neq  n-2$,
$m-1 \neq  3$.
Suppose  $m-1=3\Rightarrow m=4$. 
We wish to show there are three distinct entries:
$n-2\neq n-2m$,
$n-2 \neq  3$,
$n-2m\neq 3$.
Since $n-2m=3$ and since $m=4\Rightarrow n=11$, this  is a contradiction.
Suppose  $m-1=n-2\Rightarrow m=n-1$ and since $n-2m=3 \Rightarrow m=-2$, this  is a contradiction.
One of the following two sets has three distinct elements:  $\{m-1,n-2,3\}$ or $\{n-2,3,n-2m\}$.

\vspace{.1in}

\noindent \textit{Consider $G_2$} with  $2\leq m \leq  \frac{n}{2}-1$ and \textit{consider $G_3$} with $\frac{n+2}{6}\leq m \leq \frac{n}{2}-1$.

{\bf 1.} $G_2$ with  $2\leq m \leq  \frac{n}{2}-2$ and  $G_3$ with $\frac{n+2}{6}\leq m \leq \frac{n}{2}-2$, $l=1$, $i=2$, $j=n$.

$$\begin{pmatrix}
                    3& m-1\cr
                    n-2&n-2(m+1)\end{pmatrix}$$

We wish to show there are three distinct entries:
$n-2\neq 3$,
$m-1 \neq  n-2$,
$m-1 \neq  3$.
Suppose  $m-1=3\Rightarrow m=4$. 
We wish to show there are three distinct entries:
$n-2\neq n-2(m+1)$,
$n-2 \neq  3$,
$n-2(m+1)\neq 3$
Suppose  $n-2(m+1)=3$ and since $m=4\Rightarrow  n=13$,  which is a contradiction.
One of the following two sets has three distinct elements:    
$\{m-1,n-2,3\}$ or $\{n-2,3,n-2(m+1)\}$.

{\bf 2.} $G_2$ and $G_3$ with $m=\frac{n}{2}-1$, $l=1$, $i=2$, $j=n$
Three distinct elements  are $\{3, n-2, \frac{n}{2}-1\}$

\vspace{.1in}

\noindent \textit{Consider $G_1$,$G_2$ and $G_3$ }

 {\bf 1.} $G(i,j;l,m)$ with  $ \frac{n}{2} \leq m \leq n-2$, $l=1$, $i=2$, $j=n$.
Three distinct elements  are $\{3, n-2, n\}$.

{\bf 2.} $G(i,j;l,m)$ with $m=n-1$, $l=1$, $i=2$, $j=n$.
Three distinct elements  are $\{3, n-1, n-2\}$.

{\bf 3.} $G(i,j;l,m)$ with $m=n$, $l=1$, $i=2$, $j=n$.
Three distinct elements  are $\{3,n-2, 1\}$.
\vspace{.1in}

\noindent \emph{\textbf{Step 6}}:

For matrix $G(i,j;l,m)$ with $3 \leq i < j \leq \n-1$,  $l=1$, $2\leq m \leq n$ we look at two  cases and consider every matrix type defined above of even $n\pmod 6$.
\vspace{.1in}

\emph{\textbf{Case 1}}: We take the submatrix of $G(i,j;l,m)$ with  $l=1$, $2\leq m\leq n$, $3 \leq i \leq \frac{n}{2}+1$, $i \leq j \leq n-1$.

$$\begin{pmatrix}
                    n& \sigma^{s}{(m-1)}\cr
                    X& \sigma^{t}{(m-1)}\end{pmatrix}$$

There are three distinct entries $\{n, \sigma^{s}{(m-1)},\sigma^{t}{(m-1)}\}$.

\vspace{.1 in}

\emph{\textbf{Case 2}}: We take matrix $G(i,j;l,m)$ with $l =1$, $2 \leq m\leq n$, $\frac{n}{2}+2 \leq i <j \leq n-1$, and let $s=(n+1)-i$, $t=(n+1)-j$.

\vspace{.1in}

\noindent \textit{Consider $G_1$} 

$$\begin{pmatrix}
                    2(i-1)-n& \sigma^{s}{(m-1)}\cr
                    2(j-1)-n& \sigma^{t}{(m-1)}\end{pmatrix}$$

We wish to show there are three distinct entries:
$\sigma^{s}{(m-1)}\neq \sigma^{t}{(m-1)}$,
$2i-2-n \neq  \sigma^{s}{(m-1)}$,
$2i-2-n \neq  \sigma^{t}{(m-1)}$.
Suppose $2i-2-n \equiv  \sigma^{s}{(m-1)}\Rightarrow 3i \equiv 2n+m+2$ $\pmod {n-1}$.
Suppose $2i-2-n \equiv  \sigma^{t}{(m-1)}\Rightarrow 2i +j \equiv 2n+m+2$ $\pmod {n-1}$.
We wish to show there are three distinct entries:
$\sigma^{s}{(m-1)}\neq \sigma^{t}{(m-1)}$,
$2j-2-n \neq  \sigma^{t}{(m-1)}$,
$2j-2-n \neq  \sigma^{s}{(m-1)}$.
Suppose $2j-2-n \equiv  \sigma^{t}{(m-1)}\Rightarrow 3j \equiv 2n+m+2$.
Subtracting this equation  from $3i \equiv 2n+m+2 \Rightarrow j\equiv i$, which is a contradiction.
Subtracting the same equation from $2i +j \equiv 2n+m+2\Rightarrow j\equiv i$, which is a contradiction.
Suppose $2j-2-n \equiv  \sigma^{s}{(m-1)}\Rightarrow 2j+i \equiv 2n+m+2$.
Subtracting this equation from $3i \equiv 2n+m+2\Rightarrow  j\equiv i$, which is a contradiction.
Subtracting the same equation from $2i+j \equiv 2n+m+2\Rightarrow  j\equiv i$, which is a contradiction.
One of the following two sets has three distinct elements: 
 $\{\sigma^{s}{(m-1)},\sigma^{t}{(m-1)},2i-2-n \}$ or $\{\sigma^{s}{(m-1)},\sigma^{t}{(m-1)},2j-2-n\}$.
\vspace{.1in}

\noindent \textit{Consider $G_2$} with $l =1$, $2\leq m\leq n$, $\frac{n}{2}+2 \leq i <j \leq n-1$

{\bf 1.} $G_2$ with $i=\frac{n}{2}+2$,  $\frac{n}{2}+3 \leq j \leq n-1$, and let $s=\frac{n}{2}-1$, $t=(n+1)-j$.

$$\begin{pmatrix}
                   Y& \sigma^{s}{(m-1)}\cr
                  2(j-2)-n& \sigma^{t}{(m-1)}\end{pmatrix}$$

If  \textit{$K$ is Even} $ \Rightarrow Y=\frac{n}{2}-2$.  
We wish to show there are three distinct entries:
$\sigma^{s}{(m-1)}\neq \sigma^{t}{(m-1)}$,
$ \frac{n}{2}-2 \neq  \sigma^{t}{(m-1)}$ and
$ \frac{n}{2}-2\neq  \sigma^{s}{(m-1)}$.
Suppose $ \frac{n}{2}-2 \equiv  \sigma^{s}{(m-1)}\Rightarrow m \equiv 0$  $\pmod {n-1}$, which is a contradiction.
Suppose  $ \frac{n}{2}-2\equiv  \sigma^{t}{(m-1)}\Rightarrow 2j \equiv \n+2m+4$ $\pmod {n-1}$ and we wish to show there are three distinct entries:
$ 2(j-2)-n \neq \frac{n}{2}-2$,
$\frac{n}{2}-2 \neq  \sigma^{s}{(m-1)}$,
$2(j-2)-n\neq \sigma^{s}{(m-1)}$. 
Suppose $\frac{n}{2}-2 \equiv 2(j-2)-n \Rightarrow 4j \equiv 3n+4$. For this type  $n=6+6k\Rightarrow j\equiv 4.5k+5.5$, and since $k$ is even $j$ is not integer, which is a contradiction.
Suppose $ 2(j-2)-n \equiv  \sigma^{s}{(m-1)}\Rightarrow 4j \equiv 3n+2m+4$ $\pmod {n-1}$. 
Subtracting it from $2j \equiv \n+2m+4 \Rightarrow  j\equiv n$ $\pmod {n-1}$, which is a contradiction.
One of the following two sets has three distinct elements:  
$\{\sigma^{s}{(m-1)},\sigma^{t}{(m-1)},\frac{n}{2}-2\}$ or $\{\sigma^{s}{(m-1)},\frac{n}{2}-2,2j-4-n \}$.
\vspace{.1in}

If  \textit{$K$ is Odd} $ \Rightarrow Y=\frac{n}{2}+1$.
We wish to show there are three distinct entries:
$Y \neq \sigma^{s}{(m-1)}$,
$Y\neq  2j-n-4 $,
$2j-n-4 \neq  \sigma^{t}{(m-1)}$.
Suppose  $\frac{n}{2}+1 = 2j-n-4 \Rightarrow 4j = 3n+10 \Rightarrow j=4.5k+7$. Since $k$ is odd, $j$ is not integer, which is a contradiction.
Suppose $\frac{n}{2}+1 \equiv  \sigma^{t}{(m-1)}\Rightarrow 2j\equiv n +2m -2$ $\pmod {n-1}$.
Suppose  $ 2j-n-4  \equiv  \sigma^{t}{(m-1)}\Rightarrow 3j\equiv 2n+4+m$ $\pmod {n-1}$ and we wish to show there are three distinct entries:
$\frac{n}{2}+1 \neq 2j-n-4$,
$\frac{n}{2}+1 \neq  \sigma^{s}{(m-1)}$,
$2j-n-4 \neq  \sigma^{s}{(m-1)}$.
Suppose  $\frac{n}{2}+1 \equiv  \sigma^{s}{(m-1)}\Rightarrow m\equiv 3$ $\pmod {n-1}$. 
Subtracting this equation from  $2j\equiv n+2m-2\Rightarrow  j\equiv \frac{n}{2}+2$ $\pmod {n-1}$ and since $\frac{n}{2}+3 \leq j$ this is a contradiction.
Replacing   $m\equiv 3$ in $3j\equiv 2n+4+m \Rightarrow 3j\equiv 12k+22$, since right side is not divisible by $3$ this is a contradiction.
Suppose  $2j-n-4 \equiv  \sigma^{s}{(m-1)}\Rightarrow  4j \equiv 3n+2m+4$ $\pmod {n-1}$.
Subtracting it from $3j\equiv n+m+4\Rightarrow  j\equiv n+m\Rightarrow j>n$, which is a contradiction.
Subtracting it from $2j\equiv n+2m-2\Rightarrow j\equiv n+2m+6 \Rightarrow j>n$, which is a contradiction
One of the following two sets has three distinct elements:  
$\{\sigma^{s}{(m-1)},Y,2j-4-n \}$ or $\{Y,\sigma^{t}{(m-1)},2j-4-n \}$.

{\bf 2.} $G_2$ with $\frac{n}{2}+3 \leq i< j \leq n-1$ $l=1$, $2\leq m\leq n$, and let  $s=(n+1)-i$, $t=(n+1)-j$.

$$\begin{pmatrix}
                   2(i-2)-n& \sigma^{s}{(m-1)}\cr
                   2(j-2)-n& \sigma^{t}{(m-1)}\end{pmatrix}$$

We wish to show there are three distinct entries:
$\sigma^{s}{(m-1)}\neq \sigma^{t}{(m-1)}$,
$2i-4-n \neq  \sigma^{s}{(m-1)}$,
$2i-4-n \neq  \sigma^{t}{(m-1)}$.
Suppose $2i-4-n \equiv  \sigma^{s}{(m-1)}\Rightarrow 3i \equiv 2n+m+4$ $\pmod {n-1}$.
Suppose $2i-4-n \equiv  \sigma^{t}{(m-1)}\Rightarrow 2i +j \equiv 2n+m+4$ $\pmod {n-1}$.
We wish to show there are three distinct entries:
$\sigma^{s}{(m-1)}\neq \sigma^{t}{(m-1)}$,
$2j-4-n \neq  \sigma^{t}{(m-1)}$,
$2j-4-n \neq  \sigma^{s}{(m-1)}$.
Suppose $2j-4-n \equiv  \sigma^{t}{(m-1)}\Rightarrow 3j \equiv 2n+m+4$ $\pmod {n-1}$.
Subtracting this equation from $3i \equiv 2n+m+4 \Rightarrow  j \equiv i$, which is a contradiction.
Subtracting the same equation from $2i +j \equiv 2n+m+4 \Rightarrow j \equiv i$, which is a contradiction.
Suppose $2j-4-n \equiv  \sigma^{s}{(m-1)}\Rightarrow 2j +i \equiv 2n+m+4$ $\pmod {n-1}$.
Subtracting this equation from  $3i \equiv 2n+m+4 \Rightarrow j\equiv i$, which is a contradiction.
Subtracting the same equation from $2i+j \equiv 2n+m+4 \Rightarrow j \equiv i$, which is a contradiction.
One of the following two sets has three distinct elements: 
$\{\sigma^{s}{(m-1)},\sigma^{t}{(m-1)},2i-4-n \}$ or $\{\sigma^{s}{(m-1)},\sigma^{t}{(m-1)},2j-4-n\}$.

\vspace{.1in}

\noindent \textit{Consider $G_3$} with $l =1$, $2\leq m\leq n$, $\frac{n}{2}+2 \leq i <j \leq n-1$

{\bf  1.}  $G_3$ with  $i=\frac{n}{2}+2$,  $\frac{n}{2}+3 \leq j \leq \frac{5n+4}{6}$, $l=1$, $2 \leq m\leq n$, and let $s=\frac{n}{2}-1$,  $t=(n+1)-j$. 

$$\begin{pmatrix}
                   n-9& \sigma^{s}{(m-1)}\cr
                  2j-4-n& \sigma^{t}{(m-1)}\end{pmatrix}$$
                              
We wish to show there are three distinct entries:
$\sigma^{s}{(m-1)}\neq \sigma^{t}{(m-1)}$,
$2j-4-n \neq  \sigma^{t}{(m-1)}$,
$2j-4-n \neq  \sigma^{s}{(m-1)}$.
Suppose $2j-4-n \equiv  \sigma^{s}{(m-1)}\Rightarrow 4j \equiv 3n+2m+4$ $\pmod {n-1}$.
Suppose $2j-4-n \equiv  \sigma^{t}{(m-1)}\Rightarrow 3j\equiv 2n+m+4$ $\pmod {n-1}$. 
We wish to show there are three distinct entries:
$\sigma^{s}{(m-1)}\neq \sigma^{t}{(m-1)}$,
$n-9 \neq \sigma^{s}{(m-1)}$,
$n-9\neq  \sigma^{t}{(m-1)}$.
Suppose $n-9 \equiv  \sigma^{s}{(m-1)}\Rightarrow n \equiv  2m+14$ $\pmod {n-1}$.
Subtracting this equation from  $4j  \equiv 3n+2m+4  \Rightarrow  j\equiv\frac{n}{2}+m+4.5 \Rightarrow j$ is not integer, which is a contradiction.
Subtracting the same equation  from  $3j\equiv 2n+m+4\Rightarrow  n \equiv 3j-3m-18$. Since $n=4+6k \Rightarrow j=2k+m+7(\frac{1}{3}) \Rightarrow j $  is not integer, which is a contradiction. 
Suppose $n-9 \equiv  \sigma^{t}{(m-1)}\Rightarrow j \equiv m+9$ $\pmod {n-1}$.
Subtracting this equation from $4j \equiv 3n+2m+4\Rightarrow j \equiv  n+\frac{m-5}{3}\Rightarrow j > n$, which is a contradiction.
Subtracting the same equation  from $3j\equiv 2n+m+4 \Rightarrow  n \equiv j+2.5$, which is a contradiction.
One of the following two sets has three distinct elements: 
$\{\sigma^{s}{(m-1)},\sigma^{t}{(m-1)},2j-4-n \}$ or $\{\sigma^{s}{(m-1)},\sigma^{t}{(m-1)},n-9\}$.

Exceptions:
For $n=16$ following the argument above we have three distinct elements:
 $\{\sigma^{s}{(m-1)},\sigma^{t}{(m-1)},2j-4-n\}$ or $\{\sigma^{s}{(m-1)},\sigma^{t}{(m-1)},n-11\}$.

For $n=22$

$$\begin{pmatrix}
                   n-5& \sigma^{s}{(m-1)}\cr
                  2j-4-n& \sigma^{t}{(m-1)}\end{pmatrix}$$  

\noindent We wish to show there are three distinct entries:
$\sigma^{s}{(m-1)}\neq \sigma^{t}{(m-1)}$,
$2j-4-n \neq  \sigma^{t}{(m-1)}$,
$2j-4-n \neq  \sigma^{t}{(m-1)}$.
Suppose $2j-26\equiv  \sigma^{10}{(m-1)}\Rightarrow 2j\equiv 35+m$ $\pmod {n-1}$.
Suppose $2j-26 \equiv  \sigma^{t}{(m-1)}\Rightarrow 3j\equiv 48+m$ $\pmod {n-1}$. 
We wish to show there are three distinct entries:
$\sigma^{10}{(m-1)}\neq \sigma^{t}{(m-1)}$,
$17 \neq \sigma^{10}{(m-1)}$,
$17\neq  \sigma^{t}{(m-1)}$.
Suppose $17\equiv  \sigma^{10}{(m-1)}\Rightarrow m \equiv 8$ $\pmod {n-1}$.
Suppose $17 \equiv  \sigma^{t}{(m-1)}\Rightarrow j\equiv 5+m$ $\pmod {n-1}$. 
Subtracting two equations $m \equiv 8$ and $2j\equiv 35+m\Rightarrow j \equiv  21.5$, which is a contradiction.
Subtracting  two equations $m \equiv 8$ and $3j\equiv 48+m \Rightarrow j \equiv 17\frac{2}{3}$, which is a contradiction.
Subtracting  two equations $j\equiv 5+m$ and $2j\equiv 35+m\Rightarrow j \equiv  30$ and since $n = 22$, this is a contradiction.
Subtracting  two equations $j\equiv 5+m$ and $3j\equiv 48+m \Rightarrow j \equiv 21.5$, which is a contradiction.
One of the following two sets has three distinct elements: 
$\{\sigma^{s}{(m-1)},\sigma^{t}{(m-1)},2j-4-n\}$ or $\{\sigma^{s}{(m-1)},\sigma^{t}{(m-1)},n-5\}$.

{\bf  2.} $G_3$ with $i=\frac{n}{2}+2$, $\frac{5n+10}{6} \leq j \leq n-1$, $l=1$, $2 \leq m\leq n$, and let $t=(n+1)-j$, $s=\frac{n}{2}-1$

$$\begin{pmatrix}
                   n-9& \sigma^{s}{(m-1)}\cr
                  2j-2-n& \sigma^{t}{(m-1)}\end{pmatrix}$$
                  
We wish to show there are three distinct entries:
$\sigma^{s}{(m-1)} \neq \sigma^{t}{(m-1)}$,
$n-9 \neq  \sigma^{s}{(m-1)}$,
$n-9 \neq  \sigma^{t}{(m-1)}$.
Suppose $n-9 \equiv  \sigma^{s}{(m-1)}\Rightarrow \frac{n}{2}-7\equiv m \Rightarrow n\equiv 2m+14$  $\pmod {n-1}$.
Suppose $n-9 \equiv  \sigma^{t}{(m-1)}\Rightarrow j \equiv m+9$ $\pmod {n-1}$. 
We wish to show there are three distinct entries:
$\sigma^{s}{(m-1)}\neq \sigma^{t}{(m-1)}$,
$2j-2-n \neq  \sigma^{t}{(m-1)}$,
$2j-2-n \neq  \sigma^{s}{(m-1)}$.
Suppose $2j-2-n \equiv  \sigma^{s}{(m-1)}\Rightarrow 4j \equiv 3n+2m$ $\pmod {n-1}$.
Replacing $2m \equiv n-14 \Rightarrow  j\equiv n-3.5$, which is a contradiction.
Replacing  $m \equiv j-9$ in  $4j \equiv 3n+2m \Rightarrow j \equiv  1.5n-9 \Rightarrow j>n$, which is a contradiction. 
Suppose $2j-2-n \equiv  \sigma^{t}{(m-1)}\Rightarrow 3j\equiv 2n+2+m$ $\pmod {n-1}$. 
Replacing $m\equiv \frac{n}{2}-7$ in the equation $3j\equiv 2n+2+m \Rightarrow j\equiv \frac{5n-10}{6}$ and since $j \geq \frac{5n+10}{6}$ this is a contradiction.
Subtracting the equation $3j\equiv 2n+2+m$ from  $j \equiv m+9\Rightarrow  n \equiv j-3.5$, which is a contradiction.
One of the following two sets has three distinct elements:
$\{\sigma^{s}{(m-1)},\sigma^{t}{(m-1)},2j-2-n\}$ or $\{\sigma^{s}{(m-1)},\sigma^{t}{(m-1)},n-9\}$

Exceptions:
For $n=16$ following the argument above we have three distinct elements:
$\{\sigma^{s}{(m-1)},\sigma^{t}{(m-1)},2j-2-n\}$ or $\{\sigma^{s}{(m-1)},\sigma^{t}{(m-1)},n-11\}$.
For $n=22$ following the argument above we have three distinct elements:
$\{\sigma^{s}{(m-1)},\sigma^{t}{(m-1)},2j-2-n\}$ or $\{\sigma^{s}{(m-1)},\sigma^{t}{(m-1)},n-5\}$.

{\bf 3.} $G_3$ with $\frac{n}{2}+3 \leq i<j \leq \frac{5n+4}{6}$, $l=1$, $2 \leq m\leq n$, and let $t=(n+1)-j$,  $s=(n+1)-i$.

$$\begin{pmatrix}
                     2(i-2)-n& \sigma^{s}{(m-1)}\cr
                    2(j-2)-n& \sigma^{t}{(m-1)}\end{pmatrix}$$

This case is identical to  $G_2$ {\bf 2.}.
One of the following two sets has three distinct elements: 
$\{\sigma^{s}{(m-1)},\sigma^{t}{(m-1)},2i-4-n \}$ or $\{\sigma^{s}{(m-1)},\sigma^{t}{(m-1)},2j-4-n\}$.

{\bf 4.} $G_3$ with  $\frac{5n+10}{6}< i<j \leq n-1$ , $l=1$, $2 \leq m\leq n$, and  let $s=(n+1)-i$, $t=(n+1)-j$.

$$\begin{pmatrix}
                    2(i-1)-n& \sigma^{s}{(m-1)}\cr
                    2(j-1)-n& \sigma^{t}{(m-1)}\end{pmatrix}$$
This case is identical to  $G_1$.             
One of the following two sets has three distinct elements: 
$\{\sigma^{s}{(m-1)},\sigma^{t}{(m-1)},2i-2-n \}$ or $\{\sigma^{s}{(m-1)},\sigma^{t}{(m-1)},2j-2-n\}$.

{\bf 5.} $G_3$ with  $\frac{n}{2}+3 \leq i \leq \frac{5n+4}{6}$, $ \frac{5n+10}{6} \leq  j \leq n-1$, $l=1$, $2 \leq m\leq n$, and let $s=(n+1)-i$, $t=(n+1)-j$

$$\begin{pmatrix}
                     2(i-2)-n& \sigma^{s}{(m-1)}\cr
                    2(j-1)-n& \sigma^{t}{(m-1)}\end{pmatrix}$$

We wish to show three are  three distinct entries:
$2(i-2)-n \neq 2(j-1)-n$,
$2(i-2)-n \neq \sigma^{s}{(m-1)}$,
$2(j-1)-n \neq  \sigma^{s}{(m-1)}$.
Suppose $2i-4-n \equiv  \sigma^{s}{(m-1)}\Rightarrow 3i \equiv 2n+m+4$ $\pmod {n-1}$.
Suppose $2j-2-n \equiv  \sigma^{s}{(m-1)}\Rightarrow 2j +i \equiv 2n+m+2$ $\pmod {n-1}$.
We wish to show three are  three distinct entries:
$2i-4-n \neq 2j-2-n$,
$2j-2-n \neq  \sigma^{t}{(m-1)}$,
$2i-4-n \neq  \sigma^{t}{(m-1)}$.
Suppose $2j-2-n \equiv  \sigma^{t}{(m-1)}\Rightarrow 3j \equiv 2n+m+2$ $\pmod {n-1}$.
Subtracting this equation from $3i \equiv 2n+m+4\Rightarrow j\equiv i- \frac{2}{3}$, which is a contradiction.
Subtracting the same equation from $2i+j \equiv 2n+m+2\Rightarrow j\equiv i$, which is a contradiction.
Suppose $2i-4-n \equiv  \sigma^{t}{(m-1)}\Rightarrow j +2i \equiv 2n+m+4$ $\pmod {n-1}$.
Subtracting this equation from $3i \equiv 2n+m+4 \Rightarrow j\equiv i$, which is a contradiction.
Subtracting the same equation from $2j+i\equiv 2n+m+2\Rightarrow i \equiv j+2$, which is a contradiction.
One of the following two sets has three distinct elements:
$\{\sigma^{s}{(m-1)},2(i-2)-n,2(j-1)-n\}$ or $\{\sigma^{t}{(m-1)},2(i-2)-n, 2(j-1)-n\}$.

\vspace{.1 in}

\noindent \emph{\textbf{Step 7}}:

For matrix $G(i,j;l,m)$ with $l =1$, $2\leq m\leq n$, $3 \leq i \leq \frac{n}{2}+1$, $j=n$, and let $s=n+1-i$ and  we look at four cases and consider every matrix type defined above of even $n\pmod 6$.
\vspace{.1in}

\emph{\textbf{Case 1}}: We take matrix  $G(i,j;l,m)$ with  $l =1$, $2\leq m\leq \frac{n}{2}-1$, $3 \leq i  \leq \frac{n}{2}+1$, $j=n$,  and let $s=n+1-i$.
\vspace{.1in}

\noindent \textit{Consider  $G_1$} with $2\leq m \leq \frac{n}{2}-1$ and \textit{ $G_3$} with $2\leq m \leq \frac{n-4}{6}$.

$$\begin{pmatrix}
                    n& \sigma^{s}{(m-1)}\cr
                    n-2&n-2m\end{pmatrix}$$
               
There are three distinct elements $\{n, n-2, n-2m\}$

\vspace{.1in}

\noindent \textit{Consider  $G_2$} with $2\leq m \leq \frac{n}{2}-1$ and \textit{$G_3$} with $\frac{n+2}{6}\leq m \leq \frac{n}{2}-1$.

{\bf  1.} $G_2$  with $2\leq m \leq \frac{n}{2}-2$ and $G_3$  with $\frac{n+2}{6}\leq m \leq \frac{n}{2}-2$.

$$\begin{pmatrix}
                   n& \sigma^{s}{(m-1)}\cr
                    n-2&n-2(m+1)\end{pmatrix}$$

There are  three distinct elements $\{n, n-2, n-2(m+1)\}$

{\bf 2.} $G_2$ and $G_3$  with $m= \frac{n}{2}-1$.

$$\begin{pmatrix}
                    n& \sigma^{s}{(m-1)}\cr
                    n-2&Y\end{pmatrix}$$

For $G_2$: 
If  \textit{$K$ is Even} $ \Rightarrow Y=\frac{n}{2}-2$, there are three distinct entries $(n, n-2, \frac{n}{2}-2)$. 
If  \textit{$K$ is Odd} $ \Rightarrow Y=\frac{n}{2}+1$, there are three distinct entries $(n, n-2, \frac{n}{2}+1)$.

For $G_3$:  There are three distinct entries $\{n, n-2, n-9\}$.
Exceptions:
For $n=16$ three distinct elements are $\{16, 14, 5\}$.
For $n=22$ then there are three distinct entries $\{22, 20, 17\}$.

\vspace{.1 in}

\emph{\textbf{Case 2}}: We take matrix  $G(i,j;l,m)$ with  $l=1$, $j=n $, $3 \leq i \leq \frac{n}{2}+1$, $\frac{n}{2}\leq m \leq n-2$, and let $s=n+1-i$.

$$\begin{pmatrix}
                    n& \sigma^{s}{(m-1)}\cr
                    n-2& n\end{pmatrix}$$

Suppose $\sigma^{s}{(m-1)} \equiv n-2\Rightarrow m+2\equiv i$, since the conditions for $m$ and $i$ this is a contradiction.
There are three distinct elements $\{n,n-2,\sigma^{s}{(m-1)}\}$.
\vspace{.1 in}

\emph{\textbf{Case 3}}:  We take matrix $G(i,j;l,m)$ with  $3 \leq i \leq \frac{n}{2}+1$, $l=1$, $j=n $, $m=n-1$, and let $s=n+1-i$.

$$\begin{pmatrix}
                    n&X\cr
                    n-2& n-1\end{pmatrix}$$

Three distinct elements are $\{n,n-1,n-2\}$.
\vspace{.1 in}

\emph{\textbf{Case 4}}: We take matrix  $G(i,j;l,m)$ with  $l=1$, $j=n $, $3 \leq i \leq \frac{n}{2}+1$, $m=n$, and let $s=n+1-i$.

$$\begin{pmatrix}
                    n& X\cr
                    n-2& 1\end{pmatrix}$$

Three distinct elements are $\{n,n-2,1\}$.
\vspace{.1 in}

\noindent \emph{\textbf{Step 8}}:

For matrix $G(i,j;l,m)$ with $l =1$, $2\leq m\leq n$,  $\frac{n}{2}+2 \leq i \leq n-1$, $j=n$, and let $s=n+1-i$ and  we look at four cases and consider every matrix type defined above of even $n\pmod 6$.
\vspace{.1 in}

\emph{\textbf{Case 1}}: We take matrix $G(i,j;l,m)$ with $l =1$, $2\leq m \leq \frac{n}{2}-1$, $\frac{n}{2}+2 \leq i \leq n-1$, $j=n$, and let $s=(n+1)-i, t=(n+1)-j$.

\vspace{.1 in}

\noindent \textit{Consider  $G_1$}.

$$\begin{pmatrix}
                    2(i-1)-n& \sigma^{s}{(m-1)}\cr
                    n-2& n-2m\end{pmatrix}$$

We wish to show there are three distinct entries:
$n-2 \neq n-2m$,
$2i-2-n \neq  n-2$,
$2i-2-n  \neq  n-2m$.
Suppose $2i-2-n = n-2m\Rightarrow i=n+1+m$, which is a contradiction.
Three distinct elements are $\{n-2,n-2m, 2i-2-n\}$.

\vspace{.1in}

\noindent \textit{Consider  $G_2$}.

{\bf 1.} $G_2$ with  $2\leq m \leq \frac{n}{2}-2$,  $i=\frac{n}{2}+2$.  

$$\begin{pmatrix}
                    Y& \sigma^{s}{(m-1)}\cr
                    n-2& n-2(m+1)\end{pmatrix}$$

If  \textit{$K$ is Even} $ \Rightarrow Y=\frac{n}{2}-2 \Rightarrow$, three distinct entries $\{n-2, n-2(m+1), Y\}$.
If  \textit{$K$ is Odd} $ \Rightarrow Y=\frac{n}{2}+1 \Rightarrow$,  three distinct entries $\{n-2, n-2(m+1), Y\}$.

{\bf 2.} $G_2$ with  $2\leq m \leq \frac{n}{2}-2$, $ \frac{n}{2}+2\leq i \leq n-1$.

$$\begin{pmatrix}
                     2(i-2)-n& \sigma^{s}{(m-1)}\cr
                     n-2& n-2(m+1)\end{pmatrix}$$
                     
We wish to show there are three distinct entries:
$n-2 \neq n-2m-2$,
$2i-4-n \neq  n-2$,
$2i-4-n  \neq  n-2m-2$.
Suppose  $2i-4-n = n-2\Rightarrow i=n+1$, which is a contradiction.
Suppose  $2i-4-n = n-2-2m\Rightarrow i=n-m+2$ and since $m \leq \frac{n}{2}-2$ this is a contradiction.
Three distinct elements are  $\{n-2,n-2m-2,2i-4-n\}$.

{\bf 3.} $G_2$ with  $m=\frac{n}{2}-1$, $ \frac{n}{2}+2\leq i \leq n-1$.

$$\begin{pmatrix}
                     2(i-2)-n& \sigma^{s}{(m-1)}\cr
                     n-2& Y\end{pmatrix}$$
                  
If  \textit{$K$ is Even} $ \Rightarrow Y=\frac{n}{2}-2$.
We wish to show there are  three distinct entries:
$2i-n-4 \neq n-2$,
$n-2 \neq \frac{n}{2}-2$,
$2i-n-4 \neq \frac{n}{2}-2$.
Suppose $2i-n-4 = n-2\Rightarrow i=n+1$, which is a contradiction.
Suppose $n-2=\frac{n}{2}-2\Rightarrow n=0$, which is a contradiction.
Suppose $2i-n-4=\frac{n}{2}-2 \Rightarrow 4i-4=3n$.
We wish to show there are  three distinct entries:
$2i-n-4 \neq n-2$,
$2i-4-n \neq  \sigma^{s}{(m-1)}$,
$n-2 \neq  \sigma^{s}{(m-1)}$.
Suppose $2i-4-n \equiv  \sigma^{s}{(m-1)}\Rightarrow 6i \equiv 5n+6$ $\pmod {n-1}$.
Subtracting $4i-4=3n$  from this equation $\Rightarrow i \equiv n+1$  $\pmod {n-1}$, which is a contradiction.
Suppose $n-2 \equiv  \sigma^{s}{(m-1)}\Rightarrow 2i \equiv n+2$ $\pmod {n-1}$.
Subtracting $4i-4=3n\Rightarrow i \equiv n+1$ $\pmod {n-1}$, which is a contradiction.
One of the following two sets has three distinct elements:  
 $\{2i-n-4,n-2,Y\}$ or $\{n-2,\sigma^{s}{(m-1)},2i-4-n\}$.

If  \textit{$K$ is Odd} $ \Rightarrow Y=\frac{n}{2}+1$.
We wish to show there are  three distinct entries:
$2i-n-4 \neq n-2$,
$n-2 \neq \frac{n}{2}+1$,
$2i-n-4 \neq \frac{n}{2}+1$.
Suppose $2i-n-4 = n-2\Rightarrow i=n+1$, which is a contradiction.
Suppose  $n-2=\frac{n}{2}+1\Rightarrow n=6$, which is a contradiction.
Suppose $2i-n-4=\frac{n}{2}+1\Rightarrow 4i-10=3n$ and since $n=6k+6 \Rightarrow 18k=4i-28 \Rightarrow i=4.5k+7$ and since $k$ is odd $i$ is not integer, which is a contradiction.
There are three distinct elements $\{2i-n-4,n-2,Y\}$.

{\bf 4.} $G_2$ with $m=\frac{n}{2}-1$, $i=\frac{n}{2}+2$. 

$$\begin{pmatrix}
                    Y& \sigma^{s}{(m-1)}\cr
                    n-2& Y\end{pmatrix}$$

If  \textit{$K$ is Even} $\Rightarrow Y=\frac{n}{2}-2$.
Three distinct elements  are $\{Y,n-2,\sigma^{s}{(m-1)}\}$.
If  \textit{$K$ is Odd} $ \Rightarrow Y=\frac{n}{2}+1$. 
We wish to show there are three  distinct entries:
$Y \neq n-2$,
$Y \neq  \sigma^{s}{(m-1)}$,
$n-2 \neq  \sigma^{s}{(m-1)}$.
Suppose $Y \equiv \sigma^{s}{(m-1)}\Rightarrow 3n \equiv 8$, which is a contradiction.
Suppose  $n-2\equiv \sigma^{s}{(m-1)}\Rightarrow n\equiv m$ $\pmod {n-1}$,  which is a contradiction.
Three distinct elements $\{Y,n-2,\sigma^{s}{(m-1)}\}$.
\vspace{.1in}

\noindent \textit{Consider $G_3$}.

{\bf 1}  $G_3$ with  $i=\frac{n}{2}+2$, $j=n$,  $l=1$, $2\leq m \leq \frac{n-4}{6}$

$$\begin{pmatrix}
                    n-9& \sigma^{s}{(m-1)}\cr
                    n-2&n-2m \end{pmatrix}$$

There are three distinct entires $\{n-9, n-2, n-2m\}$.
Exceptions:
For $n=16$ three distinct entires are $\{5, 14, 16-2m\}$.
For $n=22$ three distinct entires are $\{17, 20, 22-2m\}$.

{\bf 2}  $G_3$ with  $\frac{n}{2}+3 \leq i \leq \frac{5n+4}{6}$, $j=n$, $l=1$, $2\leq m \leq \frac{n-4}{6}$. 

$$\begin{pmatrix}
                  2(i-2)-n&\sigma^{s}{(m-1)}\cr
                    n-2& n-2m\end{pmatrix}$$

We wish to show there are three distinct entries:
$n-2 \neq n-2m$,
$n-2m \neq \sigma^{s}{(m-1)}$,
$n-2 \neq \sigma^{s}{(m-1)}$.
Suppose $n-2m\equiv \sigma^{s}{(m-1)}\Rightarrow 3m \equiv i$ $\pmod {n-1}$ and this is a contradiction since $i-3m \geq 5$.
Suppose $n-2=\equiv \sigma^{s}{(m-1)}\Rightarrow i \equiv m+2$, which is a contradiction.
There are three distinct elements $\{n-2,n-2m, \sigma^{s}{(m-1)}\}$.

{\bf 3} $G_3$ with  $\frac{5n+10}{6}\leq i\leq n-1$, $j=n$, $l=1$, $2\leq m \leq \frac{n-4}{6}$. 
 
$$\begin{pmatrix}
                     2(i-1)-n& \sigma^{s}{(m-1)}\cr
                    n-2& n-2m\end{pmatrix}$$

We wish to show there are  three distinct entries:
$n-2 \neq n-2m$,
$n-2m \neq \sigma^{s}{(m-1)}$,
$n-2 \neq \sigma^{s}{(m-1)}$.
Suppose  $n-2=\equiv \sigma^{s}{(m-1)}\Rightarrow i \equiv m+2$ $\pmod {n-1}$, which is a contradiction.
Suppose $n-2m\equiv \sigma^{s}{(m-1)}\Rightarrow 3m \equiv i$ $\pmod {n-1}$ and this is a contradiction since $i-3m \geq \frac{n-1}{3}$.
Three distinct elements  $\{n-2,n-2m, \sigma^{s}{(m-1)}\}$.

{\bf 4}  $G_3$ with $l=1$, $\frac{n+2}{6}\leq m \leq \frac{n}{2}-2$, $i= \frac{n}{2}+2$, $j=n$.

$$\begin{pmatrix}
                    n-9& \sigma^{s}{(m-1)}\cr
                    n-2& n-2(m+1)\end{pmatrix}$$

Three distinct entires are $\{n-9, n-2, n-2(m+1)\}$.
Exceptions: For $n=16$ three distinct entires are $\{5, 14, 14-2m\}$.
For $n=22$ three distinct entires are $\{17, 20, 20-2m\}$.

{\bf 5}  $G_3$ with $l=1$, $\frac{n+2}{6}\leq m \leq \frac{n}{2}-2$, $\frac{n}{2}+3 \leq i \leq \frac{5n+4}{6}$, $j=n$   

$$\begin{pmatrix}
                    2(i-2)-n& \sigma^{s}{(m-1)}\cr
                    n-2& n-2(m+1)\end{pmatrix}$$

We wish to show there are  three distinct entries:
$n-2m-2 \neq n-2$,
$n-2 \neq 2i-4-n$,
$n-2m-2\neq 2i-4-n$.
Suppose $n-2= 2i-4-n\Rightarrow i=n+1$, which is a contradiction.
Suppose $n-2m-2= 2i-4-n\Rightarrow m=n-i+1$  
and we wish to show there are  three distinct entries:
$2i-4-n \neq \sigma^{s}{(m-1)} $,
$n-2\neq  2i-4-n $,
$n-2  \neq  \sigma^{s}{(m-1)}$.
Suppose  $2i-4-n =n-2\Rightarrow i=n+1$, which is a contradiction.
Suppose $n-2\equiv \sigma^{s}{(m-1)}\Rightarrow i\equiv m+2$ $\pmod {n-1}$ and since $i-m \geq 5$ this is a contradiction.
Suppose $2i-4-n \equiv  \sigma^{s}{(m-1)}\Rightarrow m\equiv 3i-2n-4$ $\pmod {n-1}$.
Subtracting this equation from $ m=n-i+1\Rightarrow 3i-2n-4\equiv n-i+1\Rightarrow 4i\equiv 3n-5$ and since $n=6k+4\Rightarrow i\equiv 4.5k+1.75 \Rightarrow i$ is not integer, which is a contradiction. 
One of the following two sets has three distinct elements: 
$\{n-2,n-2m-2,2i-4-n\}$ or $\{n-2,\sigma^{s}{(m-1)},2i-4-n\}$.

{\bf 6} $G_3$ with $l=1$, $\frac{n+2}{6}\leq m \leq \frac{n}{2}-2$, $\frac{5n+10}{6} \leq i \leq n-1$, $j=n$.

$$\begin{pmatrix}
                    2(i-1)-n&\sigma^{s}{(m-1)}\cr
                    n-2& n-2m-2\end{pmatrix}$$

We wish to show there are three distinct entries:
$n-2 \neq n-2m-2$,
$n-2m-2 \neq \sigma^{s}{(m-1)}$,
$n-2 \neq \sigma^{s}{(m-1)}$.
Suppose $n-2=\equiv \sigma^{s}{(m-1)}\Rightarrow i \equiv m+2$  $\pmod {n-1}$ and since $i-m\geq 4$ this  is a contradiction.
Suppose $n-2m-2 \equiv \sigma^{s}{(m-1)}\Rightarrow 3m \equiv i-2$ $\pmod {n-1}$.
We wish to show there are three distinct entries:
$n-2 \neq n-2m-2$,
$2i-4-n \neq  n-2$,
$2i-4-n  \neq  n-2m-2$.
Suppose  $2i-4-n = n-2m-2\Rightarrow m=n+1-i$. 
Subtracting this equation from $m \equiv \frac{i-2}{3}\Rightarrow i \equiv \frac{3n+5}{4}$ $\pmod {n-1}$and since  $i\geq \frac{5n+10}{6}$ this is a contradiction.
One of the following two sets has three distinct elements:   
$\{n-2,n-2m-2, \sigma^{s}{(m-1)}\}$ or $\{n-2,n-2m-2,2i-4-n\}$.

{\bf 7} $G_3$ with $l=1$, $m=\frac{n}{2}-1$, $i=\frac{n}{2}+2$, $j=n$

$$\begin{pmatrix}
                    n-9& \sigma^{s}{(m-1)}\cr
                    n-2& n-9\end{pmatrix}$$

Since $\sigma^{s}{(m-1)}\equiv n-3$ $\pmod {n-1}$ three distinct elements are  $(n-2, n-9,  \sigma^{s}{(m-1)})$
Exceptions: 
For $n=16$ three distinct entires are $\{14, 5, 13\}$.
For $n=22$ three distinct entires are $\{20, 17, 19\}$.

{\bf 8} $G_3$ with $m=\frac{n}{2}-1$, $\frac{n}{2}+3\leq i\leq \frac{5n+4}{6}$, $l=1$

$$\begin{pmatrix}
                    2(i-2)-n& \sigma^{s}{(m-1)}\cr
                    n-2& n-9\end{pmatrix}$$

Three distinct elements are $\{2i-n-4, n-2, n-9\}$.
Exceptions:
For $n=16$ three distinct entires are $\{2i-20, 14, 5\}$.
For $n=22$ three distinct entires are $\{2i-26, 20, 17\}$.

{\bf 9} $G_3$ with $m=\frac{n}{2}-1$, $\frac{5n+10}{6} \leq i\leq n-1$.

$$\begin{pmatrix}
                    2(i-\frac{n}{2}-1)& \sigma^{s}{(m-1)}\cr
                    n-2& n-9\end{pmatrix}$$
                
Three distinct elements are $\{2i-n-2, n-2, n-9\}$.
Exceptions:
For $n=16$ three distinct entires are $\{2i-18, 14, 5\}$.
For $n=22$ three distinct entires are $\{2i-24, 20, 17\}$.
\vspace{.1in}

\emph{\textbf{Case 2}}: $G(i,j;l,m)$ with  $\frac{n}{2}+2 \leq i \leq n-1$, $\frac{n}{2} \leq m \leq n-2$, and let $s=n+1-i$.

$$\begin{pmatrix}
                    X& \sigma^{s}{(m-1)}\cr
                    n-2& n\end{pmatrix}$$

\noindent \textit{Consider $G_1$} with  $\frac{n}{2}+2 \leq i \leq n-1$, $\frac{n}{2} \leq m \leq n-2$, $X=2(i-1)-n$.
Three distinct entries are $\{n, n-2, 2i-n-2\}$.
\vspace{.1in}

\noindent \textit{Consider $G_2$}.
 
{\bf 1.} $G_2$ with $i=\frac{n}{2}+2$.
If  \textit{$K$ is Even} $ \Rightarrow X=\frac{n}{2}-2$, there are three distinct entries $\{n, n-2, X\}$. 
If  \textit{$K$ is Odd} $ \Rightarrow X=\frac{n}{2}-1$, there are three distinct entries $\{n, n-2, X\}$.

{\bf 2.}  $G_2$ with $i \geq \frac{n}{2}+3\Rightarrow X=2(i-2)-n$.
Three  distinct entries are  $\{n, n-2, X\}$
\vspace{.1in}

\noindent \textit{Consider $G_3$}. 

{\bf 1.} $G_3$ with $i=\frac{n}{2}+2$.
There are three distinct entries  $\{n, n-2,n-9\}$

{\bf 2.} $G_3$ with $\frac{n}{2}+3\leq i\leq \frac{5n+4}{6}$.
Three distinct entries $\{n, n-2,2i-n-4\}$

{\bf 3.} $G_3$ with $\frac{5n+10}{6}\leq i\leq n-1$.
This case is similar to $G_2$.
Three  distinct entries  are $\{n, n-2, 2i-n-4\}$.
\vspace{.1in}

\emph{\textbf{Case 3}}: $G(i,j;l,m)$ with  $\frac{n}{2}+2 \leq i \leq n-1$, $j=n$, $l=1$, $m = n-1$.

$$\begin{pmatrix}
                    X&\sigma^{s}{(m-1)}\cr
                    n-2& n-1\end{pmatrix}$$

Suppose $\sigma^{n+1-i}{(n-2)}\equiv n-2\Rightarrow i+1\equiv 2n$, which is a contradiction. 
Suppose $\sigma^{n+1-i}{(n-2)}\equiv n-1\Rightarrow n\equiv i$, which is a contradiction.
Three distinct elements are $\{n-1,n-2, \sigma^{s}{(m-1)}\}$.

\vspace{.1in}

\emph{\textbf{Case 4}}: $G(i,j;l,m)$ with  $\frac{n}{2}+2 \leq i \leq n-1$, $j=n$, $l=1$, $m = n$.

$$\begin{pmatrix}
                    X& \sigma^{s}{(m-1)}\cr
                    n-2& 1\end{pmatrix}$$

Suppose $\sigma^{n+1-i}{(n-1)}\equiv n-2\Rightarrow i\equiv n+2$, which is a contradiction.
Suppose $\sigma^{n+1-i}{(n-1)}\equiv 1\Rightarrow i\equiv 2n-1$,  which is a contradiction.
Three distinct elements are $\{1,n-2, \sigma^{s}{(m-1)}\}$.

 \qed\end{proof}
\end{document}